\pdfoutput=1
\newif\ifdraft \draftfalse
\newif\ifnoname \nonamefalse
\draftfalse
\nonamefalse
\def\aistats{0} 
\ifnum\aistats=1
\draftfalse 
\documentclass[twoside,letterpaper]{article} \usepackage{aistats2017}
\else
\documentclass[11pt]{article}
\fi

\usepackage{paralist}
\usepackage{algorithm}
\usepackage{algpseudocode}
\algtext*{EndIf}
\algtext*{EndProcedure}
\ifnum\aistats=0
\usepackage{fullpage}
\fi
\usepackage{amsmath, amssymb, amsthm}
\usepackage{verbatim}
\usepackage{color}
\usepackage{bbm} 
\usepackage{xcolor}
\definecolor{DarkGreen}{rgb}{0.1,0.5,0.1}
\definecolor{DarkRed}{rgb}{0.5,0.1,0.1}
\definecolor{DarkBlue}{rgb}{0.1,0.1,0.5}
\usepackage[]{hyperref}
\hypersetup{
    unicode=false,          
    pdftoolbar=true,        
    pdfmenubar=true,        
    pdffitwindow=false,      
    pdftitle={},    
    pdfauthor={}
    pdfsubject={},   
    pdfnewwindow=true,      
    pdfkeywords={keywords}, 
    colorlinks=true,       
    linkcolor=DarkRed,          
    citecolor=DarkGreen,        
    filecolor=DarkRed,      
    urlcolor=DarkBlue,          
}
\usepackage{cleveref}
\usepackage[nottoc,notlot,notlof]{tocbibind}

\usepackage[square,numbers]{natbib}
\allowdisplaybreaks
\usepackage{tikz}

\usepackage{caption}
\usepackage{subcaption}

\newcommand\R{\mathbb{R}}

\newcommand{\cM}{\mathcal{M}}

\newcommand{\cO}{\mathcal{O}}

\newcommand{\cQ}{\mathcal{Q}}
\newcommand{\cR}{\mathcal{R}}

\newcommand{\cT}{\mathcal{T}}

\newcommand{\cX}{\mathcal{X}}

\DeclareMathOperator*{\myargmin}{\arg\!\min}

\newcommand{\bp}{\mathbf{p}}

\newcommand{\bbx}{\mathbf{x}}

\newcommand{\diag}{\ensuremath{\mathtt{Diag}}}

\newcommand{\one}{\pmb{1}}
\newcommand{\defn}{\mathrm{defn}}
\newcommand{\mult}{\mathrm{Multinomial}}

\newcommand{\GOF}{\texttt{zCDP-GOF}}
\newcommand{\MIN}{\texttt{zCDP-Min-$\chi^2$}}

\newcommand{\rendiv}[3]{\text{D}_{#1}\left(#2||#3 \right)}

\renewcommand{\hat}{\widetilde}

\DeclareMathOperator*{\Expectation}{\mathbb{E}}
\newcommand{\Ex}[2]{\Expectation_{#1}\left[#2\right]}

\newcommand{\prob}[1]{\mathrm{Pr}\left[#1\right]}

\DeclareMathOperator{\chistat}{T_\text{chi}}

\renewcommand{\hat}{\widehat}
\renewcommand{\bar}{\overline}

\newcommand{\MG}{\cM_{\text{Gauss}}}

\newcommand{\XX}[1]{X^{(#1)}}
\newcommand{\pXX}[1]{\tilde{X}^{(#1)}}

\newcommand{\UU}[1]{U^{(#1)}}
\newcommand{\pUU}[1]{U^{(n)}_{#1}}
\newcommand{\piUU}[2]{U^{(n)}_{#1,#2}}

\newcommand{\VV}[1]{V^{(#1)}}
\newcommand{\QuadChi}[1]{D^{(#1)}}
\newcommand{\estQuadChi}[1]{\hat{D}^{(#1)}}
\newcommand{\QQ}[1]{Q^{(#1)}}
\newcommand{\hQQ}[1]{\pmb{\cQ}^{(#1)}}
\newcommand{\RR}[1]{R^{(#1)}}
\newcommand{\hRR}[1]{\pmb{\cR}^{(#1)}}
\newcommand{\pipi}[1]{\pi^{(#1)}}
\newcommand{\Proj}{\pmb{P}}
\newcommand{\mintheta}{\hat\theta^{(n)}}
\newcommand{\truetheta}{\theta^0}
\newcommand{\Covar}{\Sigma}
\newcommand{\PrivCovar}[1]{\Covar_{#1}}
\newcommand{\nullbp}{\bp^0}
\newcommand{\nullp}{p^0}
\newcommand{\stat}{\mathrm{T}}

\ifnum\aistats=1
\newcommand{\mathmath}{$ }
\else
\newcommand{\mathmath}{$$}
\fi

\newtheorem{theorem}{Theorem}[section]
\newtheorem{lemma}[theorem]{Lemma}
\newtheorem{assumpt}[theorem]{Assumption}

\theoremstyle{definition}
\newtheorem{definition}[theorem]{Definition}

\theoremstyle{remark}
\newtheorem{remark}[theorem]{Remark}

\ifnum\aistats=0
\title{A New Class of Private Chi-Square Tests}
\ifnoname
\author{ Anonymous Author }
\else
\author{Daniel Kifer\thanks{Department of Computer Science and Engineering, The Pennsylvania State University.  Email: \texttt{dkifer@cse.psu.edu}.  Supported in part by NSF grant CNS-1228669.} \and Ryan Rogers\thanks{Department of Applied Mathematics and Computational Science, University of Pennsylvania. Email: \texttt{ryrogers@sas.upenn.edu}.}}
\fi
\fi
\begin{document}

\ifnum\aistats=1
%

%
\twocolumn[
\aistatstitle{A New Class of Private Chi-Square Tests}
\aistatsauthor{ Anonymous Author }
\aistatsaddress{ Unknown Institution } ]
\else
\maketitle
\fi

\begin{abstract}
In this paper, we develop new test statistics for private hypothesis testing. These statistics are designed specifically so that their asymptotic distributions, after accounting for noise added for privacy concerns, match the asymptotics of the classical (non-private) chi-square tests for testing if the multinomial data parameters lie in lower dimensional manifolds (examples include goodness of fit and independence testing). Empirically, these new test statistics outperform prior work, which focused on noisy versions of existing statistics.
%
\end{abstract}
\ifnum\aistats=0
\newpage
\tableofcontents
\newpage
\fi

\section{Introduction}
In 2008, \citet{Homer08} published a proof-of-concept attack showing that participation of individuals in scientific studies can be inferred from aggregate data typically published in genome-wide association studies (GWAS). Since then, there has been renewed interest in protecting confidentiality of participants in scientific data \cite{JS13,USF13,YFSU14,SSB16} using privacy definitions such as differential privacy and its variations \cite{DMNS06,DKMMN06,BS16,DR16}.

An important tool in statistical inference is \emph{hypothesis testing}, a general framework for determining whether a given model -- called the null hypothesis $H_0$ -- of a population should be rejected based on a sample from the population.  One of the main benefits of hypothesis testing is that it gives a way to control the probability of false discovery or Type I error -- falsely concluding that a model should be rejected when it is indeed true. 
Type II error is the probability of failing to reject $H_0$ when it is false. Typically, scientists want a test that guarantees a pre-specified Type I error (say 0.05) and has high \emph{power} -- complement of Type II error.

The standard approach to hypothesis testing is to
\begin{inparaenum}[(1)]
\item estimate the model parameters from the data, 
\item compute a \emph{test statistic} $\stat$ (a function of the data and the model parameters),
\item determine the (asymptotic) distribution of $\stat$ under the assumption that the model generated the data,
\item compute the \emph{$p$-value} (Type I error) as the probability of $\stat$ being more extreme than the realized value computed from the data.\footnote{For one-sided tests, the $p$-value is the probability of seeing the computed statistic or anything larger under $H_0$.}
\end{inparaenum}

Our main contribution is a general template for creating test statistics involving categorical data. Empirically, they improve on the power of previous work on differentially private hypothesis testing \cite{GLRV16,WLK15}, while maintaining at most some given Type I error. Our approach is to select certain properties of non-private hypothesis tests (e.g., their asymptotic distributions) and then build new test statistics that match these properties when Gaussian noise is added (e.g., to achieve \emph{concentrated differential privacy} \cite{DR16,BS16} or \emph{(approximate) differential privacy} \cite{DKMMN06}). Although the test statistics are designed with Gaussian noise in mind, other noise distributions can be applied, e.g. Laplace.\footnote{If we use Laplace noise instead, we cannot match properties like the asymptotic distribution of the non-private statistics, but the new test statistics still empirically improve the power of the tests. 
\ifnum\aistats=1
Due to space issues, these experiments appear in the supplementary materials.
\fi}

We point out that implications of this work extend beyond simply alleviating privacy concerns. In \emph{adaptive data analysis}, data may be reused for multiple analyses, each of which may depend on previous outcomes thus potentially overfitting.  This problem was recently studied in the computer science literature by Dwork et al. \cite{DFHPRR15}, who show that differential privacy can help prevent overfitting despite reusing data.  There have been several follow up works \cite{DFHPRR15Nips,CLNRW16,BNSSSU16} that improve and extend the connection between differential privacy and generalization guarantees in adaptive data analysis.  Specifically, \cite{RRST16} deals with \emph{post-selection hypothesis testing} where they can ensure a bound on Type I error even for several adaptively chosen tests, as long as each test is differentially private.


We discuss related work in Section \ref{sec:related}, provide background information about privacy in Section \ref{sec:privacy}, present our extension of minimum chi-square theory in Section \ref{sec:minchi} and show how it can be applied to goodness of fit (Section \ref{sec:gof}) and independence testing (Section \ref{sec:indep}). Experiments appear in these latter two sections.
\ifnum\aistats=0
We evaluate our test statistics with non-Gaussian noise in Section \ref{sec:lapind}.
\fi
We present conclusions in Section \ref{sec:conclusions}.

\section{Related Work}\label{sec:related}

One of the first works to study the asymptotic distributions of statistics that use differentially private data came from \citet{WZ10}.  Smith \cite{Smith11} then showed that for a large family of statistics, there is a corresponding differentially private statistic that shares the same asymptotic distribution as the original statistic.  However, these results do not ensure that statistically valid conclusions are made for finite samples. It is then the goal of a recent line of work to develop statistical inference tools that give valid conclusions for even reasonably sized datasets.  

The previous work on private statistical inference for categorical data can be roughly grouped together into two main approaches (with most primarily dealing with GWAS specific applications).  The first group adds appropriately scaled noise to the sampled data (or histogram of data) to ensure differential privacy and uses existing classical hypothesis tests, disregarding the additional noise distribution \cite{JS13}.  This approach is based on the argument that the impact of the noise becomes small as the sample size grows large.  Along these lines, \cite{VS09} studies how many more samples would be needed before the test with additional noise recovers the same level of power as the original test on the actual data.  However, as pointed out in \cite{FRY10,KS12,KS16,GLRV16}, even for moderately sized datasets, the impact of privacy noise is non-negligible and therefore such an approach can lead to misleading and statistically invalid results, specifically with much higher Type I error than the prescribed amount.  

The second group of work consists of tests that focus on adjusting step (3) in the standard approach to hypothesis testing given in the introduction.  That is, these tests use the same statistic in the classical hypothesis tests (without noise) and after making the statistic differentially private, they determine the resulting modified asymptotic distribution of the private statistic \cite{USF13,YFSU14,WLK15,GLRV16}.  Unfortunately, the resulting asymptotic distribution cannot be written analytically, and so Monte Carlo (MC) simulations or numerical approximations are commonly used to determine at what point to reject the null hypothesis.  

We focus on a different technique from these two different approaches, namely modifying step (2) in our outline of hypothesis testing.  Thus, we consider transforming the test statistic itself so that the resulting distribution is close to the original asymptotic distribution when additional Gaussian noise is used. If the noise is non-Gaussian, then this is followed by another step that appropriately adjusts the asymptotic distribution.  The idea of modifying the test statistic for  \emph{regression coefficients} to obtain a $t$-statistic in ordinary least squares has also been considered in \cite{Or15}.

\section{Privacy Preliminaries}\label{sec:privacy}
Formal privacy definitions can be used to protect scientific data with the careful injection of noise. Hypothesis testing must then properly account for this noise to avoid generating false conclusions. Since our primary focus is on the noise added by the privacy definitions (rather than their specific privacy semantics), we briefly discuss the privacy definitions and then elaborate on how to add noise to satisfy those definitions.

Let $\cX$ be an arbitrary domain for records. We define two datasets $\bbx= (x_1,\cdots, x_n), \bbx'=(x_1',\cdots, x_n') \in \cX^n$ to be \emph{neighboring} if they differ in at most one entry, i.e. there is some $i \in [n]$ where $x_i \neq x_i'$, but $x_j = x_j'$ for all $j \neq i$. We now define differential privacy (DP)\cite{DMNS06,DKMMN06}.

\begin{definition}[Differential Privacy]
A randomized algoirthm $\cM: \cX^n \to \cO$ is $(\epsilon,\delta)$-DP if for all neighboring datasets $\bbx,\bbx'$ and each subset of outcomes $S \subseteq \cO$, 
$$
\prob{\cM(\bbx) \in S} \leq e^{\epsilon}\prob{\cM(\bbx') \in S} + \delta.
$$
If $\delta = 0$, we simply say $\cM$ is $\epsilon$-DP.
\end{definition}


\ifnum\aistats=1
In this work, we focus on a recent variation of differential privacy, called \emph{zero concentrated differential privacy} (zCDP) \cite{BS16}; extensions of our work to $\epsilon$-DP can be found in the supplementary material due to lack of space. Further, due to space restrictions and complexity of zCDP, its formal definition can also be found in the full version of this paper.  The following result shows that zCDP lies between \emph{pure}-DP where $\delta =0$ and \emph{approximate}-DP where $\delta$ may be positive.
\begin{theorem}[\citep{BS16}]\label{thm:zcdpcompare}
If $\cM$ is $\epsilon$-DP, then $\cM$ is $\frac{\epsilon^2}{2}$-zCDP.  Further, if $\cM$ is $\rho$-zCDP then $\cM$ is $(\rho + 2\sqrt{\rho\ln(1/\delta)},\delta)$-DP for every $\delta>0$.  
\label{thm:reduction}
\end{theorem}

The following property is useful because it ensures the privacy of the dataset no matter what an adversary does with the output of a zCDP algorithm.
\begin{theorem}[Post Processing \citep{BS16}]\label{thm:post}
Let $\cM: \cX^n \to \cO$ and $g: \cO \to \cO'$ be randomized algorithms.  If $\cM$ is $\rho$-zCDP then $g \circ \cM: \cX^n \to \cO'$ is $\rho$-zCDP.  
\end{theorem}

We can privately release a function $f: \cX^n \to \R^d$ of the data using the Gaussian Mechanism $\MG$ \cite{DKMMN06}. $\MG$ first computes the \emph{global sensitivity} of $f$, which is defined as 
\mathmath
\Delta_p(f) = \max_{\text{neighboring } \bbx, \bbx' \in \cX^n }\{ ||f(\bbx) - f(\bbx') ||_p\}.
\mathmath
and then generates a noisy version of $f$ as follows (here $\sigma=\Delta_2(f)/\sqrt{2\rho}$):
\begin{equation}
\MG(\bbx) \sim N(f(\bbx),\sigma^2 \cdot I_d).
\label{eq:GM}
\end{equation}

\begin{theorem}[\cite{BS16}]
For a function $f: \cX^n \to \R^d$, the Gaussian mechanism $\MG$ from \eqref{eq:GM} is $\rho$-zCDP.
\label{thm:gauss_mech}
\end{theorem}
\else
In this work, we focus on a recent variation of differential privacy, called \emph{zero concentrated differential privacy} (zCDP) \cite{BS16} (see also \cite{DR16} for the definition of concentrated differential privacy which \cite{BS16} modified).  In order to define zCDP, we first define the R\'{e}nyi-divergence between two probability distributions.
\begin{definition}[R\'{e}nyi-Divergence]
Let $P_1$ and $P_2$ be probability distributions on space $\Omega$.  For $\alpha \in (1,\infty)$, we define the R\'{e}nyi-Divergence of order $\alpha$ of $P_1$ from $P_2$ as 
$$
\rendiv{\alpha}{P_1}{P_2} = \frac{1}{\alpha-1} \log\left( \Ex{x \sim P_1}{\left( \frac{P_1(x)}{P_2(x)}\right)^{\alpha-1}}\right).
$$
\end{definition}
\begin{remark}
Note as $\alpha \to 1$ we get KL-divergence and as $\alpha \to \infty$ we get max-divergence.
\end{remark}

We are now ready to define zCDP.
\begin{definition}[zCDP]
A mechanism $\cM: \cX^n \to \cO$ is $\rho$-zero concentrated differentially private (zCDP), if for all neighboring datasets $\bbx,\bbx' \in \cX^n$  and all $\alpha \in (1,\infty)$ we have 
$$
\rendiv{\alpha}{\cM(\bbx)}{\cM(\bbx')} \leq \rho \alpha
$$
\end{definition}The following result shows that zCDP lies between \emph{pure}-DP where $\delta =0$ and \emph{approximate}-DP where $\delta$ may be positive.
\begin{theorem}[\citep{BS16}]\label{thm:zcdpcompare}
If $\cM$ is $\epsilon$-DP, then $\cM$ is $\frac{\epsilon^2}{2}$-zCDP.  Further, if $\cM$ is $\rho$-zCDP then $\cM$ is $(\rho + 2\sqrt{\rho\ln(1/\delta)},\delta)$-DP for every $\delta>0$.  
\label{thm:reduction}
\end{theorem}

In order to compute some statistic $f: \cX^n \to \R^d$ on the data, a differentially private mechanism is to simply add symmetric noise to $f(\bbx)$ with standard deviation that depends on the \emph{global sensitivity} of $f$, which we define as 
$$
\Delta_p(f) = \max_{\text{neighboring } \bbx, \bbx' \in \cX^n }\{ ||f(\bbx) - f(\bbx') ||_p\}.
$$

In statistical hypothesis tests, it is typical to use the central limit theorem to form statistics of the data that are asymptotically normally distributed.  Then we can determine whether to reject the given model in hypothesis testing by computing the corresponding $p$-values based on the asymptotic distribution of the statistic, which works well in practice.  Because Gaussian random variables have nice composition guarantees, like the sum of two Gaussian random variables is again Gaussian (a property that is not shared with Laplace random variables), it is then desirable to use a privacy definition which is more accustomed to Gaussian perturbations.  We then define the Gaussian mechanism $\cM_{Gauss}: \cX^n  \to \R^d$ for statistic $f: \cX^n \to \R^d$, where $\sigma = \frac{\Delta_2(f) }{\sqrt{2\rho}}$, as
\begin{equation}
\MG(\bbx) \sim N(f(\bbx),\sigma^2I_d).
\label{eq:GM}
\end{equation}

\begin{theorem}
For statistic $f: \cX^n \to \R^d$, the Gaussian mechanism $\MG$ is $\rho$-zCDP.
\label{thm:gauss_mech}
\end{theorem}

We now state several of the nice properties that zCDP shares with DP.

\begin{theorem}[Post Processing \citep{BS16}]\label{thm:post}
Let $\cM: \cX^n \to \cO$ and $g: \cO \to \cO'$ be randomized algorithms.  If $\cM$ is $\rho$-zCDP then $\cM': \cX^n \to \cO'$ where $\cM'(\bbx) = g(\cM(\bbx))$ is $\rho$-zCDP.  
\end{theorem}
\begin{theorem}[Composition \citep{BS16}]
Let $\cM_1: \cX^n \to \cO$ and $\cM_2: \cX^n \to \cO'$ be randomized algorithms where $\cM_1$ is $\rho_1$-zCDP and $\cM_2$ is $\rho_2$-zCDP.  Then the composition $\cM: \cX^n \to \cO\times \cO'$ where $\cM(\bbx) = (\cM_1(\bbx), \cM_2(\bbx))$ is $(\rho_1+\rho_2)$-zCDP.  
\end{theorem}

\fi
For this work we will be considering categorical data. 
That is, we assume the domain $\cX$ has been partitioned into $d$ buckets or outcomes
and the function $f: \cX^n \to \R^d$ returns a histogram counting how many records are in each bucket. Our test statistics will only depend on this histogram. Since neighboring datasets $\bbx, \bbx^\prime$ of size $n$ differ on only one entry, their corresponding histograms differ by $\pm 1$ in exactly two buckets.  Hence, we will say that two histograms are neighboring if they differ in at most two entries by at most 1.
In this case, \ifnum\aistats=0$ \Delta_1(f) = 2$ and \fi $\Delta_2(f) = \sqrt{2}$.  To preserve privacy, we will add noise to the corresponding histogram $X = (X_1,\cdots, X_d)$ of our original dataset to get $\tilde{X}=(\tilde{X}_1,\dots,\tilde{X}_d)$. We perform hypothesis testing on this noisy histogram $\tilde{X}$.  By \Cref{thm:post}, we know that each of our hypothesis tests will be $\rho$-zCDP as long as we add Gaussian noise with variance  $1/\rho$ to each count in $X$.  
\ifnum\aistats=0
Similarly, if we add Laplace noise with scale $2/\epsilon$ to each count, we will achieve $\epsilon$-DP (this is just an instance of the Laplace Mechanism \cite{DMNS06}).
\fi

\section{General Chi-Square Tests}\label{sec:minchi}
In the non-private setting, a chi-square test involves a histogram $X$ and a model $H_0$ that produces expected counts $\bar{X}$ over the $d$ buckets. In general, $H_0$ will have $k < d$ parameters and will estimate the parameters from $X$. The chi-square test statistic is defined as
\mathmath
\chistat=\sum_{i=1}^d (X_i-\bar{X}_i)^2/\bar{X}_i.
\mathmath
If the data were generated from $H_0$ and if $k$ parameters had to be estimated, then the asymptotic distribution of $\chistat$ is $\chi^2_{d-k-1}$, a chi-square random variable with $d-k-1$ degrees of freedom. This is the property we want our statistics to have when they are computed from the noisy histogram $\tilde{X}$ instead of $X$.  Note that in the classical chi-square tests (e.g. Pearson independence test), the statistic $\chistat$ is computed and if it is larger than the $1-\alpha$ percentile of $\chi^2_{d-k-1}$, then the model is rejected.

The above facts are part of a more general \emph{minimum chi-square asymptotic theory} \cite{Ferg96}, which we overview in Section \ref{subsec:minchi}. However, we first explain the differences between private and non-private asymptotics \cite{WLK15,GLRV16}.

\subsection{Private Asymptotics}\label{subsec:asymptotics}
In \underline{non-private statistics}, a function of $n$ data records is considered a random variable, and non-private asymptotics considers this distribution as $n\rightarrow\infty$. 
In \underline{private asymptotics}, there is another quantity $\sigma^2_n$, the variance of the added noise. 

In the \emph{classical private regime}, one studies what happens as $n/\sigma^2_n\rightarrow\infty$; i.e., when the variance due to privacy is insignificant compared to sampling variance in the data (i.e. $O(n))$. In practice, asymptotic distributions derived under this regime result in unreliable hypothesis tests because privacy noise is significant \cite{USF13}.

In the \emph{variance-aware private regime}, one studies what happens as $n/\sigma^2_n\rightarrow$ constant as $n\rightarrow\infty$; that is, when the variance due to privacy is proportional to sampling variance. In practice, asymptotic distributions derived under this regime result in hypothesis tests with reliable Type I error (i.e. the $p$-values they generate are accurate) \cite{GLRV16,WLK15}.
From now on, we will be using the variance-aware privacy regime.\footnote{Note that taking $n$ and $\sigma_n^2$ to infinity is just a mathematical tool for simplifying expressions while mathematically keeping privacy noise variance proportional to the data variance; it does not mean that the amount of actual noise added to the data  depends on the data size.}

\subsection{Minimum Chi-Square Theory}\label{subsec:minchi}
In this section, we present important results about \emph{minimum chi-square theory}.
The discussion is based largely on \cite{Ferg96} (Chapter 23). Our work relies on this theory to construct new private test statistics in Sections \ref{sec:gof} and \ref{sec:indep} whose asymptotic behavior matches the non-private asymptotic behavior of the classical chi-square test.

We consider a sequence of $d$-dimensional random vectors $\VV{n}$ for $n \geq 1$ (e.g. the data histogram). The parameter space $\Theta$ is a non-empty open subset of $\R^k$, where $k\leq d$. The model $A$ maps a $k$-dimensional parameter $\theta\in \Theta$ into a $d$-dimensional vector (e.g., the expected value of $\VV{n}$), hence it maps $\Theta$ to a subset of a $k$-dimensional manifold in $d$-dimensional space.

In this abstract setting, the \underline{null hypothesis} is that there exists a $\truetheta\in\Theta$ such that:\footnote{Here $\stackrel{D}{\to}$ means convergence in distribution, as in the Central Limit Theorem \cite{Ferg96}.} 
\begin{equation}
\sqrt{n}\left(\VV{n} - A(\truetheta) \right) \stackrel{D}{\to} N(0, C(\truetheta))
\label{eq:asym_norm}
\end{equation}
where $C(\theta) \in \R^{d\times d}$ is a covariance matrix.
Intuitively, Equation \ref{eq:asym_norm} says that the Central Limit Theorem can be applied for $\truetheta$.

 We measure the distance between $\VV{n}$ and $A(\theta)$ with a test statistic given by the following quadratic form:
\begin{equation}
\QuadChi{n}(\theta) = n \left(\VV{n} - A(\theta) \right)^\intercal \hspace{-0.3em}M(\theta) \left( \VV{n} - A(\theta)\right)
\label{eq:quadratic}
\end{equation}
where $M(\theta) \in \R^{d\times d}$ is a symmetric positive-semidefinite matrix; different choices of $M$ will result in different test statistics.
We make the following standard assumptions about $A(\theta)$ and $M(\theta)$. 
\begin{assumpt} For all $\theta \in \Theta$, we have:
\ifnum\aistats=0
\begin{itemize}
\item \else 1) \fi $A(\theta)$ is bicontinuous,\footnote{i.e. $\theta_j\rightarrow\theta\Leftrightarrow A(\theta_j)\rightarrow A(\theta)$.}
\ifnum\aistats=0 \item \else 2) \fi  $A(\theta)$ has continuous first partial derivatives, which we denote as $\dot{ A}(\theta)$ with full rank $k$,
\ifnum\aistats=0 \item \else 3) \fi  $M(\theta)$ is continuous in $\theta$ and there exists an $\eta >0$ such that $M(\theta) - \eta I_d$ is positive definite in an open neighborhood of $\truetheta$.
\ifnum\aistats=0
\end{itemize}
\fi
\label{assumpt:mean}
\end{assumpt} 

\ifnum\aistats=0
The following theorem will be useful in determining the distribution for the quadratic form $\QuadChi{n}(\theta)$.
\begin{theorem}[\cite{Ferg96}]
Let $W \sim N(0,\Lambda)$. $W^\intercal W \sim \chi_r^2$ (chi-square distribution with $r$ degrees of freedom) if and only if $\Lambda$ is a projection of rank $r$. If $\Lambda$ is invertible, $W^\intercal\Lambda^{-1}W\sim\chi_r^2$.
\label{thm:chi2}
\end{theorem}
\fi
If $\truetheta$ is known, \ifnum\aistats=1 then we show in the supplementary file that \fi setting $M(\theta)  = C(\theta)^{-1}$ in  \eqref{eq:quadratic} \ifnum\aistats=0 and applying \Cref{thm:chi2} shows that \fi then $\QuadChi{n}(\truetheta)$ converges in distribution to $\chi^2_{d}$. However, as we show in \Cref{sec:gof}, this can be a sub-optimal choice of $M$.

When $\truetheta$ is not known, we need to estimate a good parameter $\mintheta$ to plug into \eqref{eq:quadratic}. One approach is to set $\mintheta = \arg\min_{\theta\in\Theta} \QuadChi{n}(\theta)$. However, this can be a difficult optimization.
If there is a rough estimate of $\truetheta$ based on the data, call it $\phi(\VV{n})$, and if it converges in probability to $\truetheta$ (i.e. $\phi(\VV{n}) \stackrel{P}{\to} \truetheta \text{ as } n \to \infty$), then we can plug it into the middle matrix to get:
\begin{equation}
\estQuadChi{n}(\theta) =  n \left(\VV{n} - A(\theta) \right)^\intercal M(\phi(\VV{n}) ) \left( \VV{n} - A(\theta)\right).
\label{eq:R}
\end{equation}
and then set our estimator $\mintheta=\arg\min_{\theta\in\Theta} \estQuadChi{n}(\theta)$. The test statistic becomes $\estQuadChi{n}(\mintheta)$ and the following theorems describe its asymptotic properties under the null hypothesis.
We use the shorthand $A = A(\truetheta)$, $M = M(\truetheta)$, and $C = C(\truetheta)$.  See the appendix for the full proof, which follows a similar argument as in \cite{Ferg96}.
\begin{theorem}
Let $\mintheta = \myargmin_{\theta \in \Theta} \estQuadChi{n}(\theta)$. 
Given \Cref{assumpt:mean} and \eqref{eq:asym_norm}, we have $\sqrt{n}(\mintheta - \truetheta)  \stackrel{D}{\to} N(0,\Psi)$ where $\truetheta$ is the true parameter and 
\mathmath
\Psi = \left(\dot{A}^\intercal M \dot{A} \right)^{-1} \dot{A}^\intercal M C M \dot{A} \left(\dot{A}^\intercal M \dot{A} \right)^{-1}.
\mathmath
\label{thm:min_chi}
\end{theorem}

We then state the following result using a slight modification of Theorem 24 in \cite{Ferg96}, which we prove in the \ifnum\aistats=1 supplementary file\else appendix\fi.
\begin{theorem}
Let $\nu$ be the rank of $C(\theta_0)$.
If \Cref{assumpt:mean} and \eqref{eq:asym_norm} hold, and, for all $\theta \in \Theta$,
\mathmath
C(\theta)M(\theta)C(\theta) = C(\theta)
\mathmath
and
\mathmath
 C(\theta)M(\theta)  \dot A(\theta) = \dot A(\theta)
\mathmath
then for $\mintheta$ given in \Cref{thm:min_chi} and $\estQuadChi{n}(\theta)$ given in \eqref{eq:R} we have: 
\mathmath
\estQuadChi{n}\left(\mintheta\right) \stackrel{D}{\to} \chi^2_{\nu-k}
\mathmath
\label{thm:ferg}
\end{theorem}

\section{Private Goodness of Fit Tests}\label{sec:gof}
As a warmup, we will first cover goodness of fit testing where the null hypothesis is simply testing whether the underlying unknown parameter is equal to a particular value.  We consider categorical data $\XX{n} = \left(\XX{n}_1,\cdots, \XX{n}_d\right) \sim \mult(n,\bp)$ where $\bp = (p_1,\cdots, p_d)$ is some probability vector over the $d$ outcomes.  We want to test the null hypothesis $H_0: \bp = \nullbp$, where each component of $\nullbp$ is positive, but we want to do so in a private way.  We then have the following classical result \citep{BFH75}.

\begin{lemma}
Under the null hypothesis $H_0: \bp = \nullbp$, $\XX{n}/n$ is asymptotically normal
\mathmath
\sqrt{n} \left(\frac{\XX{n}}{n} - \nullbp\right) \stackrel{D}{\to} N(0,\Covar)
\mathmath
where $\Covar$ has rank $d-1$ and can be written as 
\begin{equation}
\Covar \stackrel{\defn}{=} \diag(\nullbp) - \nullbp (\nullbp)^\intercal.
\label{eq:Sigma}
\end{equation}
\end{lemma}

\subsection{Unprojected Private Test Statistic}
To preserve $\rho$-zCDP, we will add appropriately scaled Gaussian noise to each component of the histogram $\XX{n}$.  We then define the zCDP statistic $\pUU{\rho} = \left(\piUU{\rho}{1}, \cdots, \piUU{\rho}{d} \right)$ where we write $Z \sim N\left(0,1/\rho \cdot I_d\right)$ and
\begin{equation}
\pUU{\rho} \stackrel{\defn}{=} \sqrt{n}\left( \frac{\XX{n} + Z}{n} - \nullbp\right).
\label{eq:UzCDP}
\end{equation}

We next derive \ifnum\aistats=1 (see proof in supplementary file) \fi the asymptotic distribution of $\pUU{\rho}$ under both private asymptotic regimes in Section \ref{subsec:asymptotics} (note that $\sigma^2=1/\rho)$.
\begin{lemma}\label{lem:dist_gof}
The random vector $\pUU{\rho_n}$ from \eqref{eq:UzCDP} under the null hypothesis $H_0: \bp = \nullbp$ has the following asymptotic distribution. If $n\rho_n \to \infty$ then $\pUU{\rho_n} \stackrel{D}{\to} N(0,\Covar)$.  Further, if $n \rho_n \to \rho>0$ then $\pUU{\rho_n} \stackrel{D}{\to} N(0,\PrivCovar{\rho})$ where $\PrivCovar{\rho}$  has full rank and 
\begin{equation}
\PrivCovar{\rho} \stackrel{\defn}{=} \Covar + 1/\rho \cdot I_d.
\label{eq:covar_gof}
\end{equation}
\end{lemma}
\ifnum\aistats=0
\begin{proof}
We know from the central limit theorem that $\pUU{\rho}$ will converge in distribution to a multivariate normal with covariance matrix given in \eqref{eq:covar_gof}.  We now show that $\PrivCovar{\rho}$ is full rank.  From \eqref{eq:Sigma} we know that $\Sigma$ is positive-semidefinite because it is a covariance matrix, hence it has all nonnegative eigenvalues.  We then consider the eigenvalues of $\PrivCovar{\rho}$. Let $\bbx \in \R^d$ be be an eigenvector of $\PrivCovar{\rho}$ with eigenvalue $\lambda \in \R$, i.e.
$$
\PrivCovar{\rho} \bbx = \lambda \bbx \implies  \Covar\bbx = (\lambda - 1/\rho)\bbx.
$$  
We then must have that $\bbx$ is also an eigenvector of $\Sigma$.  Because $\Sigma$ is positive-semidefinite we have the following inequality
$$
\lambda - 1/\rho \geq 0 \implies \lambda \geq 1/\rho>0.
$$
Thus, all the eigenvalues of $\PrivCovar{\rho}$ are positive, which results in $\PrivCovar{\rho}$ being nonsingular.
\end{proof}
\fi

Because  $\PrivCovar{\rho}$ is invertible when the privacy parameter $\rho> 0$, we can create a new statistic based on $\pUU{\rho}$ that has a chi-square asymptotic distribution under variance-aware privacy asymptotics.
\begin{theorem}
Let $\pUU{\rho_n}$ be given in \eqref{eq:UzCDP} for $n \rho_n \to\rho >0$.  If the null hypothesis $H_0: \bp = \nullbp$ holds, then for $\Sigma_{ n\rho_n}$ given in \eqref{eq:covar_gof}, we have
\begin{equation}
\QQ{n}_{\rho_n} \stackrel{\defn}{=} \left(\pUU{\rho_n}\right)^\intercal \PrivCovar{n \rho_n}^{-1} \pUU{\rho_n} \stackrel{D}{\to} \chi^2_d.
\label{eq:new_stat1}
\end{equation}
\end{theorem}
\ifnum\aistats=0
\begin{proof}
We directly apply \Cref{thm:chi2} with $W^{(n)} = \PrivCovar{n \rho_n}^{-1/2}\pUU{\rho_n}$ which is asymptotically multivariate normal with mean zero and covariance $\PrivCovar{\rho}^{-1/2}\PrivCovar{\rho}\PrivCovar{\rho}^{-1/2} = I_d$.  
\end{proof}
\fi

\ifnum\aistats=0
By computing the inverse of $\PrivCovar{n\rho_n}$ we can simplify the statistic  $\QQ{n}_{\rho_n}$. 

\begin{lemma}
We can rewrite the statistic in \eqref{eq:new_stat1} as
\begin{equation}
\QQ{n}_{\rho} =   \sum\limits_{i=1}^d \frac{\left(\piUU{\rho}{i}\right)^2}{\nullp_i+\frac{1}{n\rho}} +  \frac{n\rho}{\sum_{\ell=1}^d \frac{\nullp_\ell}{\nullp_\ell + \frac{1}{n\rho}} }\left(\sum\limits_{j=1}^d \frac{\nullp_j}{\nullp_j+ \frac{1}{n\rho}} \cdot\piUU{\rho}{j}\right)^2.
\label{eq:simple_stat1}
\end{equation}
\end{lemma}
\begin{proof}
We begin by writing the inverse of the covariance matrix $\PrivCovar{\rho}$ from \eqref{eq:covar_gof} by applying Woodbury's formula \citep{Wood50} which gives the inverse of a modified rank deficient matrix,
\begin{equation}
\PrivCovar{\rho}^{-1} = \diag(\nullbp+1/\rho\cdot \one)^{-1} + \frac{1}{1-\nullbp\cdot\omega(\rho)} \omega(\rho) \omega(\rho)^\intercal 
\label{eq:covar_inverse}
\end{equation}
where $\omega(\rho) =  \left( \frac{\nullp_1}{\nullp_1 + 1/\rho}, \cdots, \frac{\nullp_d}{\nullp_d + 1/\rho}\right)^\intercal = \frac{\nullbp}{\nullbp+1/\rho\cdot \one}$.

We note that the vector $\one$ is an eigenvector of $\Sigma_{\rho}$ and $\Sigma_\rho^{-1}$ with eigenvalue $1/\rho$ and $\rho$, respectively.  Letting $\pXX{n}_i = \XX{n}_i + Z_i$ be the perturbed version of $\XX{n}_i$ leads to the test statistic
\begin{align*}
\left(\pUU{\rho}\right)^\intercal \PrivCovar{n \rho}^{-1} \pUU{\rho} &=\sum\limits_{i=1}^d  \frac{(\pXX{n}_i-n\nullp_i)^2}{n\nullp_i+1/\rho} +  \frac{1}{1 - \sum_i \frac{(\nullp_i)^2}{\nullp_i+\frac{1}{n\rho}}}\left(\sum\limits_{i=1}^d \frac{(\pXX{n}_i-n\nullp_i)}{\sqrt{n}}\frac{\nullp_i}{\nullp_i+\frac{1}{n\rho}}\right)^2  \\
& = \sum\limits_{i=1}^d \frac{(\pXX{n}_i-n\nullp_i)^2}{n\nullp_i+1/\rho} +  \frac{1}{1 - \sum_i \frac{(\nullp_i)^2}{\nullp_i+\frac{1}{n \rho} }}\left(\sum\limits_{i=1}^d \frac{(\pXX{n}_i-n\nullp_i)}{\sqrt{n}}\frac{\nullp_i}{\nullp_i+\frac{1}{n \rho}}\right)^2
\end{align*}
We can then rewrite the term in the denominator,
$$ 
1- \sum_{i=1}^d\frac{(\nullp_i)^2}{\nullp_i+\frac{1}{n \rho}} =\sum_{i=1}^d\left(\frac{\nullp_i(\nullp_i+\frac{1}{n \rho})}{\nullp_i+\frac{1}{n \rho}}-\frac{(\nullp_i)^2}{\nullp_i+\frac{1}{n \rho}} \right) = \frac{1}{n \rho} \cdot \sum_{i=1}^d\frac{\nullp_i}{\nullp_i+\frac{1}{n \rho}}.
$$
Recalling the form of $\pUU{\rho}$ from \eqref{eq:UzCDP} concludes the proof.
\end{proof}

Note that the coefficient on the second term of \eqref{eq:simple_stat1} grows large as $n \rho\rightarrow \infty$, so this test statistic does not approach the nonprivate test for a fixed $\rho$. This is not surprising since $\PrivCovar{n \rho }$ must converge to a singular matrix as $n \rho\rightarrow \infty$.  
\else
Note that $\PrivCovar{\rho}$ is ill-conditioned when $\rho$ is large (data variance overwhelms privacy noise), since $\Covar$ is singular and $\Covar\ \one=\mathbf{0}$. This makes the test statistic unstable. 
\fi
  Further, the additional noise adds a degree of freedom to the asymptotic distribution of the original statistic.  This additional degree of freedom results in increasing the point in which we reject the null hypothesis, i.e. the critical value.  Thus, rejecting an incorrect model becomes harder as we increase the degrees of freedom, and hence decreases power.  

\subsection{Projected Private Test Statistic}
Given that the test statistic in the previous section depends on a nearly singular matrix,
we now derive a new test statistic for the private goodness of fit test. It has the remarkable property that its asymptotic distribution is $\chi^2_{d-1}$ under both private asymptotics. 

We start with the following observation.  In the classical chi-square test, the random variables $\left(\frac{(\XX{n}_i-n\nullp_i)}{\sqrt{n\nullp_i}}\right)_{i=1}^d$ have covariance matrix $I_d - \sqrt{\nullbp}\sqrt{\nullbp}^\intercal$ under the null hypothesis $H_0: \bp = \nullbp$. The classical test essentially uncorrelates these random variables and projects 
 them onto the subspace orthogonal to $\sqrt{\nullbp}$. We will use a similar intuition for the privacy-preserving random vector $\pUU{\rho}$.
 

The matrix $\PrivCovar{\rho}$ in  \eqref{eq:covar_gof} has eigenvector $\one$ with eigenvalue $1/\rho$ -- regardless of the true parameters of the data-generating distribution. Hence we think of this direction as pure noise.
We therefore project $\pUU{\rho}$ onto the space orthogonal to $\one$ (i.e. enforce the constraint that the entries in $\pUU{\rho}$ add up to $0$, as they would in the noiseless case).  We then define the \emph{projected statistic} $\hQQ{n}_\rho$ as the following where we write the projection matrix $\Proj \stackrel{\defn}{=} I_d - \frac{1}{d} \one\one^\intercal$
\begin{equation}
 \hQQ{n}_{\rho} \stackrel{\defn}{=} \left(\pUU{\rho}\right)^\intercal\Proj \PrivCovar{n \rho }^{-1} \Proj \pUU{\rho}.
\label{eq:new_stat2}
\end{equation}

It will be useful to write out the middle matrix in  $\hQQ{n}_{\rho_n}$ for analyzing its asymptotic distribution \ifnum\aistats=0 which we prove in the supplementary file\fi.
\begin{lemma}
For the covariance matrix $\PrivCovar{n\rho_n}$ given in \eqref{eq:covar_gof}, we have the following identity when $n \rho_n \to \rho >0$
\ifnum\aistats=1
$
\else
\begin{equation}
\fi
\Proj \PrivCovar{n \rho_n }^{-1} \Proj \to \PrivCovar{\rho}^{-1} - \frac{ \rho }{d} \cdot \one\one^\intercal 
\ifnum\aistats=1
$
\else
\label{eq:middle_rho}
\end{equation}
\fi
Further, when $n \rho_n \to \infty$, we have the following 
\ifnum\aistats=1 
$ 
\else 
\begin{equation} 
\fi
 \Proj \PrivCovar{n \rho_n }^{-1} \Proj \to \Proj\diag \left(\nullbp \right)^{-1}\Proj
\ifnum\aistats=1 
$
\else 
\label{eq:middle_infty}
\end{equation}
\fi
\end{lemma}
\ifnum\aistats=0
\begin{proof}
To prove \eqref{eq:middle_rho}, we use the fact that $\PrivCovar{\rho}^{-1}$ has eigenvalue $\rho$ with eigenvector $\one$.  We then focus on proving \eqref{eq:middle_infty}.  We use the identity for the inverse of $\PrivCovar{n\rho_n}^{-1}$ from \eqref{eq:covar_inverse}.  
\begin{align*}
& \Proj \PrivCovar{n \rho_n }^{-1} \Proj \\
&\qquad  = \Proj \diag(\nullbp + \frac{1}{n\rho_n} \cdot \one)^{-1}\Proj \\
& \qquad \qquad  + \frac{n \rho_n}{\sum_{i=1}^d \frac{\nullp_i}{\nullp_i + \frac{1}{n\rho_n}}} \cdot \Proj \left( \frac{\nullbp}{\nullbp + \frac{1}{n\rho_n}\one}\right)\left( \frac{\nullbp}{\nullbp + \frac{1}{n\rho_n}\one}\right)^\intercal \Proj
\end{align*}
We then focus on the second term and write $\lambda_n = \sum_{i=1}^d \frac{\nullp_i}{\nullp_i + \frac{1}{n\rho_n}}$.
\begin{align*}
& \frac{n \rho_n}{\lambda_n}\cdot \Proj \left( \frac{\nullbp}{\nullbp + \frac{1}{n\rho_n}\one}\right)\left( \frac{\nullbp}{\nullbp + \frac{1}{n\rho_n}\one}\right)^\intercal \Proj\\
&\qquad  = \frac{n \rho_n}{\lambda_n} \cdot  \left( \frac{\nullbp}{\nullbp + \frac{1}{n\rho_n}\one} - \frac{\lambda_n}{d}\cdot  \one \right)\left( \frac{\nullbp}{\nullbp + \frac{1}{n\rho_n}\cdot \one} - \frac{\lambda_n}{d}\cdot \one \right)^\intercal \\
& \qquad = \frac{n \rho_n}{\lambda_n} \cdot \left( \frac{\nullbp}{\nullbp + \frac{1}{n\rho_n}\cdot \one} \right)\left( \frac{\nullbp}{\nullbp + \frac{1}{n\rho_n}\cdot \one} \right)^\intercal \\
& \qquad \qquad - \frac{n\rho_n}{d}\cdot  \left( \frac{\nullbp}{\nullbp + \frac{1}{n\rho_n}\cdot \one} \right) \one^\intercal - \frac{n\rho_n}{d}\cdot \one \left( \frac{\nullbp}{\nullbp + \frac{1}{n\rho_n}\cdot \one} \right)^\intercal + \frac{n  \rho_n \lambda_n}{d^2} \one\one^\intercal
\end{align*}
We consider entry $(i,j)$ of the above matrix, which we can write as 
\begin{align*}
& \frac{n \rho_n}{\lambda_n} \cdot  \left( \frac{\nullp_i}{\nullp_i + \frac{1}{n\rho_n}} \right)\cdot \left( \frac{\nullp_j}{\nullp_j + \frac{1}{n\rho_n}} \right) - \frac{n \rho_n}{d} \left(\frac{\nullp_i}{\nullp_i+\frac{1}{n\rho_n}} + \frac{\nullp_j}{\nullp_j+\frac{1}{n\rho_n}} \right) +  \frac{n \rho_n \lambda_n}{d^2}\\
& \qquad = \frac{n \rho_n}{d \lambda_n} \left(\frac{\lambda_n^2}{d} - \frac{1}{(\nullp_i+\frac{1}{n\rho_n})(\nullp_j + \frac{1}{n \rho_n})} \left( \frac{\lambda_n}{n \rho_n}(\nullp_i +\nullp_j) - \nullp_i\nullp_j (d - 2 \lambda_n) \right) \right) \\
& \qquad = n \rho_n \left(\frac{\lambda_n}{d^2} - \frac{(2\lambda_n - d)\nullp_i\nullp_j}{d(\nullp_i+\frac{1}{n\rho_n})(\nullp_j + \frac{1}{n \rho_n})} \right) - \frac{\nullp_i + \nullp_j}{d \lambda_n(\nullp_i+\frac{1}{n\rho_n})(\nullp_j + \frac{1}{n \rho_n})}.
\end{align*}
We then let $n \to \infty$ to get
$$
\frac{n \rho_n}{d} \left(\frac{\lambda_n}{d} - \frac{(2\lambda_n - d)\nullp_i\nullp_j}{\lambda_n(\nullp_i+\frac{1}{n\rho_n})(\nullp_j + \frac{1}{n \rho_n})} \right) - \frac{\nullp_i + \nullp_j}{d \lambda_n(\nullp_i+\frac{1}{n\rho_n})(\nullp_j + \frac{1}{n \rho_n})} \to \frac{1}{\nullp_i} + \frac{1}{\nullp_j} - \frac{1}{\nullp_i} - \frac{1}{\nullp_j} = 0.
$$
Thus, we have shown that for $n \rho_n \to \infty$,
\begin{align*}
\Proj \PrivCovar{n \rho_n }^{-1} \Proj \to \Proj \diag\left( \nullbp \right)^{-1} \Proj.
\end{align*}

\end{proof}
\fi

\ifnum\aistats=0 We now show that the projected statistic is asymptotically chi-square distributed in both private asymptotic regimes.
\else We give the full proof of the following asymptotic distribution for the projected statistic in both private asymptotic regimes in the supplementary file. 
\fi

\ifnum\aistats=1
\begin{theorem}
Let $\pUU{\rho}$ be given in \eqref{eq:UzCDP}.  The projected statistic $\hQQ{n}_{\rho}$ has the following asymptotic distribution for $ n\rho_n \to \rho>0$ and $n \rho_n \to \infty$ (as $n\rightarrow\infty$) if the null hypothesis $H_0: \bp = \nullbp$ holds: 
$\hQQ{n}_{\rho_n} \stackrel{D}{\to} \chi^2_{d-1}$.    
\label{thm:proj_asympt}
\end{theorem}
\else
\begin{theorem}
Let $\pUU{\rho}$ be given in \eqref{eq:UzCDP}. For null hypothesis $H_0: \bp = \nullbp$, we can write the projected statistic $\hQQ{n}_{\rho}$ in the following way for $\tilde n = \sum_{i=1}^d (\XX{n}_i + Z_i)$
\begin{align}
\hQQ{n}_{\rho} & =   \sum\limits_{i=1}^d \frac{\left(U_i^{(n)}(\rho)\right)^2}{p_i^0+\frac{1}{n\rho}} - \frac{\rho}{d} \left( \tilde n - n\right)^2\nonumber \\ 
& + \frac{n\rho}{\sum_{\ell=1}^d  \frac{p_\ell^0}{p_\ell^0 + \frac{1}{n\rho}}}\left(\sum\limits_{j=1}^d  \frac{p_j^0}{p_j^0 + \frac{1}{n\rho}}\cdot \UU{n}_j(\rho)\right)^2 .
\label{eq:simple_stat2}
\end{align}
Further for $n\rho_n \to \rho>0$, if the null hypothesis holds then $\hQQ{n}_{\rho_n} \stackrel{D}{\to} \chi^2_{d-1}$.  
\label{thm:proj_asympt}
\end{theorem}
\begin{proof}
We first show that we can write the projected statistic in \eqref{eq:new_stat2} in the proposed way.  Using \eqref{eq:middle_rho}, we can write the projected statistic in terms of the unprojected statistic in \eqref{eq:simple_stat1}, which will give the expression in \eqref{eq:simple_stat2}
$$
\hQQ{n}_{\rho_n} =  U_i^{(n)}(\rho_n)^\intercal \left( \Sigma_{n\rho_n}^{-1} - \frac{ n\rho_n }{d} \cdot \one\one^\intercal \right)U_i^{(n)}(\rho_n) = \QQ{n}_{\rho_n} - \frac{n \rho_n }{d} \cdot U_i^{(n)}(\rho_n)^\intercal \one\one^\intercal U_i^{(n)}(\rho_n).
$$

We then turn to determining the asymptotic distribution of the projected statistics when $n \rho_n \to \rho >0$.  Recall that $\one$ is an eigenvector of $\PrivCovar{\rho}$.  Note that $\PrivCovar{\rho}$ is diagonalizable, i.e. $\PrivCovar{\rho} = BDB^\intercal$ where $D$ is a diagonal matrix and $B$ is an orthogonal matrix with one column being $1/d \cdot \one$.  For the following matrix $\Lambda$, we can write it as a $d\times d$ identity matrix except one of the entries on the diagonal is zero.
$$
\Lambda = \PrivCovar{\rho}^{-1/2} \Proj BDB^\intercal\Proj \PrivCovar{\rho}^{-1/2}.
$$
Thus, $\Lambda$ is idempotent and has rank $d-1$.  We define $W \sim N(0,I_{d-1})$.   We then know that $\hQQ{n}_{\rho_n}$ has the same asymptotic distribution as $W^\intercal W$ and so we can apply \Cref{thm:chi2}.    
\end{proof}

\begin{theorem}
For histogram data $\XX{n}$, the projected statistic $\hQQ{n}_{\rho_n}$ in \Cref{thm:proj_asympt} converges in distribution to a $\chi^2_{d-1}$ when $H_0: \bp = \nullbp$ holds and $n \rho_n \to \infty$.  In fact, as $n \rho_n \to \infty$, the difference between $\hQQ{n}_{\rho_n}$ and the classical chi-square statistic $\sum_{i=1}^d \frac{\left(\XX{n}_i - n \nullp_i\right)^2}{n\nullp_i}$ converges in probability to 0.
\label{cor:proj}
\end{theorem}
\begin{proof}
Although $\PrivCovar{n\rho_n}^{-1}$ does not exist as $n \rho_n \to \infty$, we can still write the asymptotic projected statistic.  The middle matrix in the projected statistic when $n \rho_n \to \infty$ is then $\Proj\diag(\nullbp)^{-1}\Proj$.

When $n \rho_n \to \infty$, we also have that $\pUU{\rho_n} \stackrel{D}{\to} N(0,\Covar)$ from \Cref{lem:dist_gof}.  We then analyze the asymptotic distribution of the projected statistic, where we write $U \sim N(0,\Covar)$ and study the distribution of $U^\intercal\Proj\diag(\nullbp)^{-1} \Proj U$.  We note that we have $U^\intercal \one = 0$, which simplifies the asymptotic distribution of the projected statistic.
\begin{align*}
U^\intercal\Proj\diag(\nullbp)^{-1} \Proj  U & = \sum_{i=1}^d \frac{U_i^2}{\nullp_i} 
\end{align*}
Note that this last final form is exactly the original chi-square statistic used in the classical test, which is known to converge to $\chi^2_{d-1}$. 
\end{proof}
\fi

\subsection{Comparison of Statistics}
We now want to compare the two private chi-square statistics in \eqref{eq:new_stat1} and \eqref{eq:new_stat2} to see which may lead to a larger \emph{power} (i.e. smaller Type II error).  The following theorem shows that we can write the unprojected statistic \eqref{eq:new_stat1} as a combination of both the projected statistic \eqref{eq:new_stat2} and squared independent Gaussian noise.
\begin{theorem}
Consider histogram data $\XX{n}$ that has Gaussian noise $Z \sim N(0,1/\rho\cdot I_d)$ added to it.  For the statistics $\QQ{n}_\rho$ and $\hQQ{n}_{\rho}$ based on the noisy counts given in \eqref{eq:new_stat1} and \eqref{eq:new_stat2} respectively, we have
\mathmath
\QQ{n}_\rho = \hQQ{n}_{\rho} + \frac{n\rho}{d} \left(\sum_{i=1}^d Z_i \right)^2.
\mathmath
Further, for any fixed data $\XX{n}$, $\hQQ{n}_{\rho} $ is independent of $\left(\sum_{i=1}^d Z_i \right)^2$.  
\label{thm:compare}
\end{theorem}

\ifnum\aistats=1
We present the full proof in the supplementary file.
\else
To prove this we will use the noncentral version of Craig's Theorem.
\begin{theorem}[Craig's Theorem \citep{RD88}]
Let $Y \sim N(\mu,V)$.  Then the quadratic forms $Y^\intercal AY$ and $Y^\intercal BY$ are independent if $AVB = 0$.  
\end{theorem}
We are now ready to prove our theorem.
\begin{proof}[Proof of \Cref{thm:compare}]
We first show that we can write $\QQ{n}_\rho - \hQQ{n}_{\rho} = \frac{n\rho}{d} \left(\sum_{i=1}^d Z_i \right)^2$.  Note that $\left(\pUU{\rho}\right)^\intercal\one = \sum_{i=1}^d Z_i/\sqrt{n}$ and $\PrivCovar{\rho}^{-1}$ has eigenvalue $\rho$ with eigenvector $\one$.  We then have
\begin{align*}
\QQ{n}_{\rho}& = \left(\pUU{\rho}\right)^\intercal\PrivCovar{n\rho}^{-1} \pUU{\rho} \\
& = \left(\pUU{\rho}\right)^\intercal\left(I_d - \frac{1}{d}\one\one^\intercal + \frac{1}{d}\one\one^\intercal \right)^\intercal\PrivCovar{n\rho}^{-1} \left(I_d - \frac{1}{d}\one\one^\intercal + \frac{1}{d}\one\one^\intercal \right)\pUU{\rho} \\
& = \hQQ{n}_{\rho}   + \frac{2}{d} (\pUU{\rho})^\intercal \one \one^\intercal\PrivCovar{n\rho}^{-1}\left(I_d - \frac{1}{d}\one\one^\intercal\right)\pUU{\rho} 
+ \frac{1}{d^2}(\pUU{\rho})^\intercal \one \one^\intercal\PrivCovar{n\rho}^{-1}  \one\one^\intercal \pUU{\rho} \\
& = \hQQ{n}_{\rho} + \frac{2}{d} (\pUU{\rho})^\intercal \one \one^\intercal\PrivCovar{n\rho}^{-1}\Proj\pUU{\rho} + \frac{n\rho}{d}\left( \sum_{i=1}^d Z_i\right)^2 \\
& = \hQQ{n}_{\rho} + \frac{2n \rho}{d} \left(\sum_{i=1}^d Z_i / \sqrt{n} \right) \one^\intercal\pUU{\rho} - \frac{2n \rho}{d} \left(\sum_{i=1}^d Z_i\right)^2 + \frac{n\rho}{d}\left( \sum_{i=1}^d Z_i\right)^2 \\
& = \hQQ{n}_{\rho} + \frac{n\rho}{d} \left(\sum_{i=1}^d Z_i \right)^2
\end{align*}

We now apply Craig's Theorem to show that for a fixed histogram $\XX{n}$, we have $\hQQ{n}_{\rho} $ is independent of $\left(\sum_{i=1}^d Z_i \right)^2$. When $\XX{n}$ is fixed, we can define the random variable $Y \sim N\left(\mu, 1/\rho I_d \right)$ where $\mu = (\XX{n}- n\nullbp)/\sqrt{n}$.  If we set $A = \Proj\PrivCovar{n\rho}^{-1}\Proj$, then our projected statistic can be rewritten as $Y^\intercal AY$.  Further, if we define $B = \one\one^\intercal$, then $\left(\sum_{i=1}^d Y_i\right)^2 = Y^\intercal BY$.  We then have $A\left( 1/\rho\cdot  I_d\right) B = 0$, so that the projected statistic is independent of $\left(\sum_{i=1}^d Y_i\right)^2$.  We next note that $Y = Z + \mu$ and that $\one^\intercal\mu = 0$.  Hence, 
$$
Y^\intercal BY = (Z+ \mu)^\intercal B (Z+\mu) = Z^\intercal BZ + 2\mu^\intercal B Z + \mu^\intercal B \mu = Z^\intercal B Z = \left(\sum_{i=1}^d Z_i\right)^2.
$$
\end{proof}
\fi

Algorithm \ref{alg:gof} ($\GOF$) shows how to perform goodness of fit testing with either of these two test statistics, i.e. unprojected \eqref{eq:new_stat1} or projected \eqref{eq:new_stat2}.
%
%
%
We note that our test is zCDP for neighboring histogram datasets due to it being an application of the Gaussian mechanism and \Cref{thm:post}.  Hence:
\begin{theorem}
$\GOF(\cdot;\rho,\alpha,\nullbp)$ is $\rho$-zCDP.
\end{theorem}

\begin{algorithm}
\caption{ zCDP Chi-Square Goodness of Fit Test}
\label{alg:gof}
\begin{algorithmic}
\Procedure{\GOF}{$\XX{n}$; $\rho$, $\alpha$, $H_0: \bp = \nullbp$}
\\ Set $\pXX{n}\gets \XX{n} + Z$ where $Z \sim N(0,1/\rho \cdot I_d)$.
\State \textbf{For the \emph{unprojected statistic}}:
\State $
\stat \gets \frac{1}{n}\left(\pXX{n} -n \nullbp \right)^\intercal \PrivCovar{n\rho}^{-1} \left(\pXX{n} -n \nullbp\right)
$
\State $t\gets (1-\alpha)$ quantile of $\chi^2_d$
\State \textbf{For the \emph{projected statistic}}:
\State $
\stat \gets \frac{1}{n}  \left(\pXX{n} - n\nullbp \right)^\intercal \Proj \PrivCovar{n\rho}^{-1} \Proj \left(\pXX{n} - n\nullbp \right)
$
\State $t\gets (1-\alpha)$ quantile of $\chi^2_{d-1}$\\
{\bf if }{$\stat > t$} {\bf then } 
  Reject
\EndProcedure
\end{algorithmic}
\end{algorithm}
 
\ifnum\aistats=1
When the null hypothesis is false (i.e., $\bp\neq\nullbp$), both statistics converge to a non-central chi-square distribution (the analysis can be found in the supplementary file).  We then turn to empirical results.
\else
\subsection{Power Analysis}
From \Cref{thm:compare} we see that the difference between $\QQ{n}_{\rho}$ and $\hQQ{n}_{\rho}$ is the addition of squared independent noise.  This additional noise can only hurt \emph{power}, because for the same data the statistic $\QQ{n}_{\rho}$ has larger variance than $\hQQ{n}_{\rho}$ and does not depend on the underlying data.  If we fix an alternate hypothesis, we can obtain asymptotic distributions for our two test statistics.
\begin{theorem}
Consider the null hypothesis $H_0: \bp = \nullbp$ and the alternate hypothesis $H_1:\bp = \nullbp + \frac{\pmb{\Delta}}{\sqrt{n}}$ where $\sum_{i=1}^d\Delta_i = 0$.  Assuming the data $\XX{n}$ comes from the alternate $H_1$, the two statistics $\hQQ{n}_{\rho_n}$, and $\QQ{n}_{\rho_n}$ have noncentral chi-square distributions when $n \rho_n  \to \rho > 0$, i.e. 
\mathmath
\QQ{n}_{\rho_n}\stackrel{D}{\to} \chi_{d}^2\left(\pmb{\Delta}^\intercal \PrivCovar{\rho}^{-1}\pmb{\Delta}\right) \quad \& \quad 
\hQQ{n}_{\rho_n}\stackrel{D}{\to} \chi_{d-1}^2\left(\pmb{\Delta}^\intercal \PrivCovar{\rho}^{-1}\pmb{\Delta}\right) .
\mathmath
Further, if $n \rho_n \to \infty$ then 
\mathmath
\hQQ{n}_{\rho_n}\stackrel{D}{\to} \chi_{d-1}^2\left(\sum_i \frac{\Delta_i^2}{\nullp_i} \right)
\mathmath
\label{thm:power}
\end{theorem}

We point out that in the case where $n \rho_n \to \infty$, the projected statistic has the same asymptotic distribution as the classical (nonprivate) chi-square test under the same alternate hypothesis.  

We will use the following result to prove this theorem.
\begin{lemma}[\cite{Ferg96}]
Suppose $Y  \sim N(\pmb{\delta}, V)$.  If $V$ is a projection of rank $\nu$ and $V \pmb{\delta} = \pmb{\delta}$ then $U^\intercal U \sim \chi^2_\nu(\pmb{\delta}^\intercal\pmb{\delta})$.  
\label{lem:Ferg_noncentral}
\end{lemma}

\begin{proof}[Proof of \Cref{thm:power}]
In this case we have the random vector $\pUU{\rho_n}$ from \eqref{eq:UzCDP} converging in distribution to $N(\pmb{\Delta},\PrivCovar{\rho})$ if $n \rho_n \to \rho>0$ or $N(\pmb{\Delta},\Covar)$ if $n \rho_n \to \infty$ by \Cref{lem:dist_gof}.  We first consider the case when $n \rho_n \to \rho>0$.  Consider $U \sim N(\pmb{\Delta},\PrivCovar{\rho})$ and $Y = \PrivCovar{\rho}^{-1/2} U \sim N((\PrivCovar{\rho}^{-1/2}\pmb{\Delta},I_d)$.  We then know that $Y^\intercal Y$ and the unprojected statistic $\QQ{n}_{\rho_n}$ have the same asymptotic distribution.  In order to use \Cref{lem:Ferg_noncentral}, we need to verify that $\PrivCovar{\rho}^{-1/2}\PrivCovar{\rho}\PrivCovar{\rho}^{-1/2}\left(\PrivCovar{\rho}^{-1/2}\pmb{\Delta}\right) = \PrivCovar{\rho}^{-1/2}\pmb{\Delta}$, which indeed holds.

We then consider the projected statistic $\hQQ{n}_{\rho_n}$ where $n\rho_n \to \rho>0$.  Similar to the proof of \Cref{thm:proj_asympt}, we diagonalize $\PrivCovar{\rho} = B D B^\intercal$ where $B$ is an orthogonal matrix with one column being $1/d \cdot \one$ and $D$ is a diagonal matrix.  We then let $U\sim N(\pmb{\Delta},\PrivCovar{\rho})$ and let 
$$
Y = \PrivCovar{\rho}^{-1/2} \Proj U
$$
We then have $Y^\intercal Y$ and $\hQQ{n}_{\rho_n}$ will have the same asymptotic distribution.  Recall that $\Lambda = \PrivCovar{\rho}^{-1/2} \Proj \PrivCovar{\rho} \Proj \PrivCovar{\rho}^{-1/2}$ is idempotent with rank $d-1$.  Lastly, to apply \Cref{lem:Ferg_noncentral} we need to show the following
$$
\Lambda \left(\PrivCovar{\rho}^{-1/2} \Proj \pmb{\Delta} \right)=\PrivCovar{\rho}^{-1/2} \Proj \pmb{\Delta}.
$$
Let $\hat B \in \R^{d\times (d-1)}$ be the same as matrix $B$ whose corresponding column for $1/d \cdot \one$ is missing, which we assume to be the last column of $B$.  Further, we define $\hat D \in \R^{(d-1)\times (d-1)}$ to be the same as $D$ without the last row and column.  We can then write $\Proj \Sigma_{\rho} \Proj = \hat B \hat D\hat B^\intercal$.  We can then simplify the left hand side to have 
\begin{align*}
\PrivCovar{\rho}^{-1/2}&  \Proj \PrivCovar{\rho} \Proj \PrivCovar{\rho}^{-1} \Proj \\
& = B D^{-1/2}B^\intercal \Proj \PrivCovar{\rho}\hat B \hat D^{-1} \hat B^\intercal \\
& = B D^{-1/2}B^\intercal \hat B D B^\intercal \hat B \hat D^{-1} \hat B^\intercal \\
& = B D^{-1/2}B^\intercal \hat B D \hat D^{-1} \hat B^\intercal \\
& = B D^{-1/2} \hat B^\intercal \\
& = B D^{-1/2}  B^\intercal \Proj \\
& = \PrivCovar{\rho}^{-1/2} \Proj 
\end{align*}
The noncentral parameter is then 
$$
\pmb{\Delta}^\intercal \Proj \PrivCovar{\rho}^{-1} \Proj \pmb{\Delta}
$$
We then note that $\sum_i \Delta_i = 0$.

For the case when $n \rho_n \to \infty$.  From \eqref{eq:middle_infty}, we have $\Proj \PrivCovar{n \rho_n}\Proj \to M_\infty$, which can be diagonalized.  As we showed in \Cref{cor:proj}, we have 
$$
\left(\pUU{\rho_n}\right)^\intercal M_\infty \pUU{\rho_n} =  \left(\pUU{\rho_n}\right)^\intercal \diag(\nullbp)^{-1} \pUU{\rho_n}
$$
From \Cref{lem:dist_gof}, we know that $\left(\pUU{\rho_n}\right)\stackrel{D}{\to} N(\pmb{\Delta},\Covar)$.  We then write $U \sim N(\pmb{\Delta},\Covar)$ so that our projected chi-square statistic has the same asymptotic distribution as 
$$
U^\intercal\diag(\nullbp)^{-1} U
$$
which has a $\chi^2_{d-1}(\pmb{\Delta}\diag(\nullbp)^{-1} \pmb{\Delta})$ distribution.
\end{proof}

Note that the noncentral parameters in the previous theorem are the same for both statistics and only the degrees of freedom are different.  
\fi


\subsection{Experiments for Goodness of Fit Testing}\label{sect:GOF_Results}
Throughout all of our experiments, we will fix $\alpha = 0.05$ and privacy parameter $\rho = 0.001$.  All of our tests are designed to achieve Type I error at most 
$\alpha$\ifnum\aistats=1.\footnote{Due to space limitations we give the empirical Type I error for various $\nullbp$ and $n$ in the supplementary file. }
\else, as we empirically show for different null hypotheses $\nullbp$ and sample size $n$ in \Cref{fig:GOF_signif}.  We include 1.96 times the standard error of our $100,000$ independent trials (giving a $95\%$ confidence interval) for each sample size and each null hypothesis.
\fi

\ifnum\aistats=0
\begin{figure}
\begin{center}
\begin{subfigure}{.45\textwidth}
\begin{flushleft}
\includegraphics[width=\linewidth]{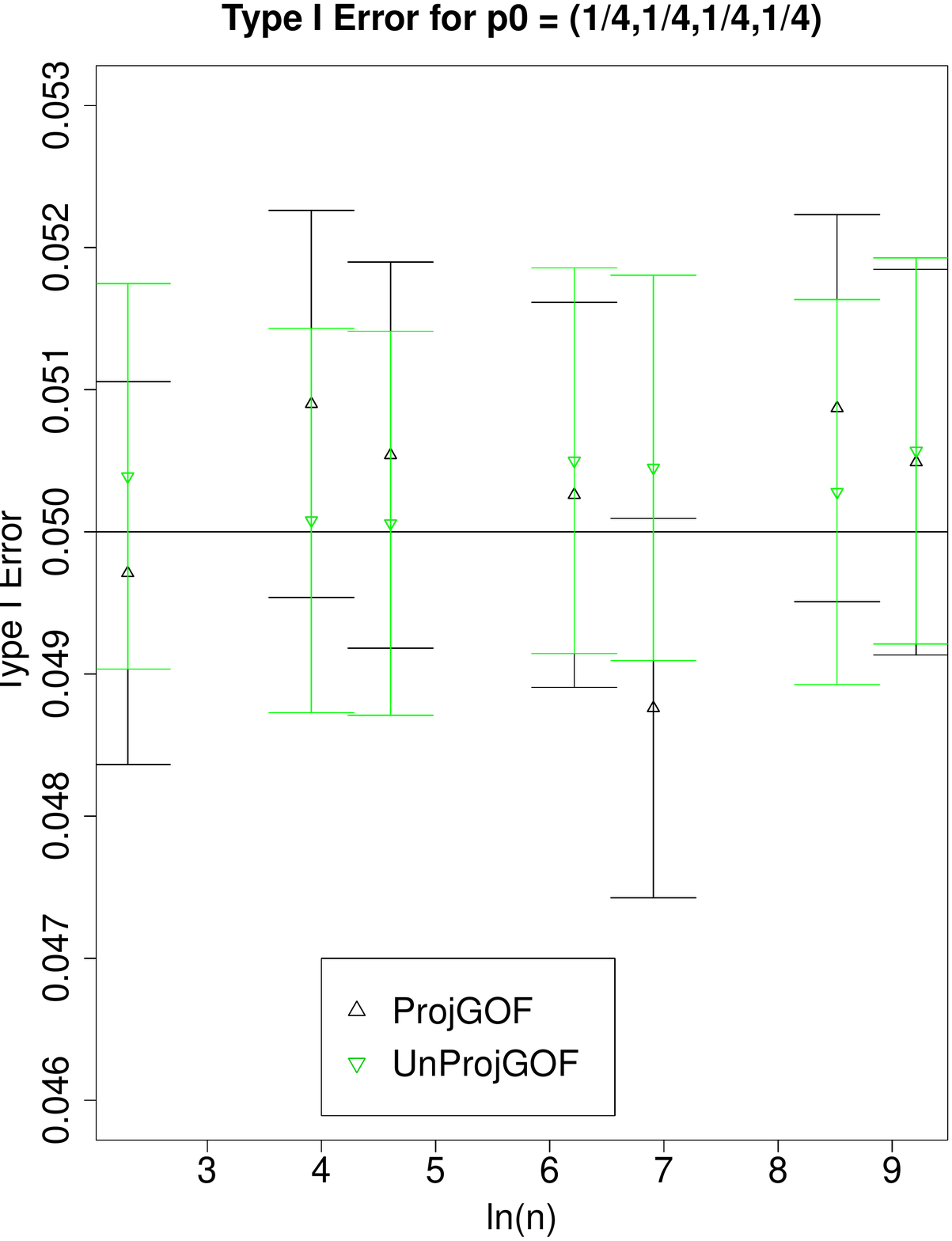}
\end{flushleft}
\end{subfigure}
\hspace{10mm}
\begin{subfigure}{.1\textwidth}
\end{subfigure}
\begin{subfigure}{.45\textwidth}
\begin{flushright}
\includegraphics[width=\linewidth]{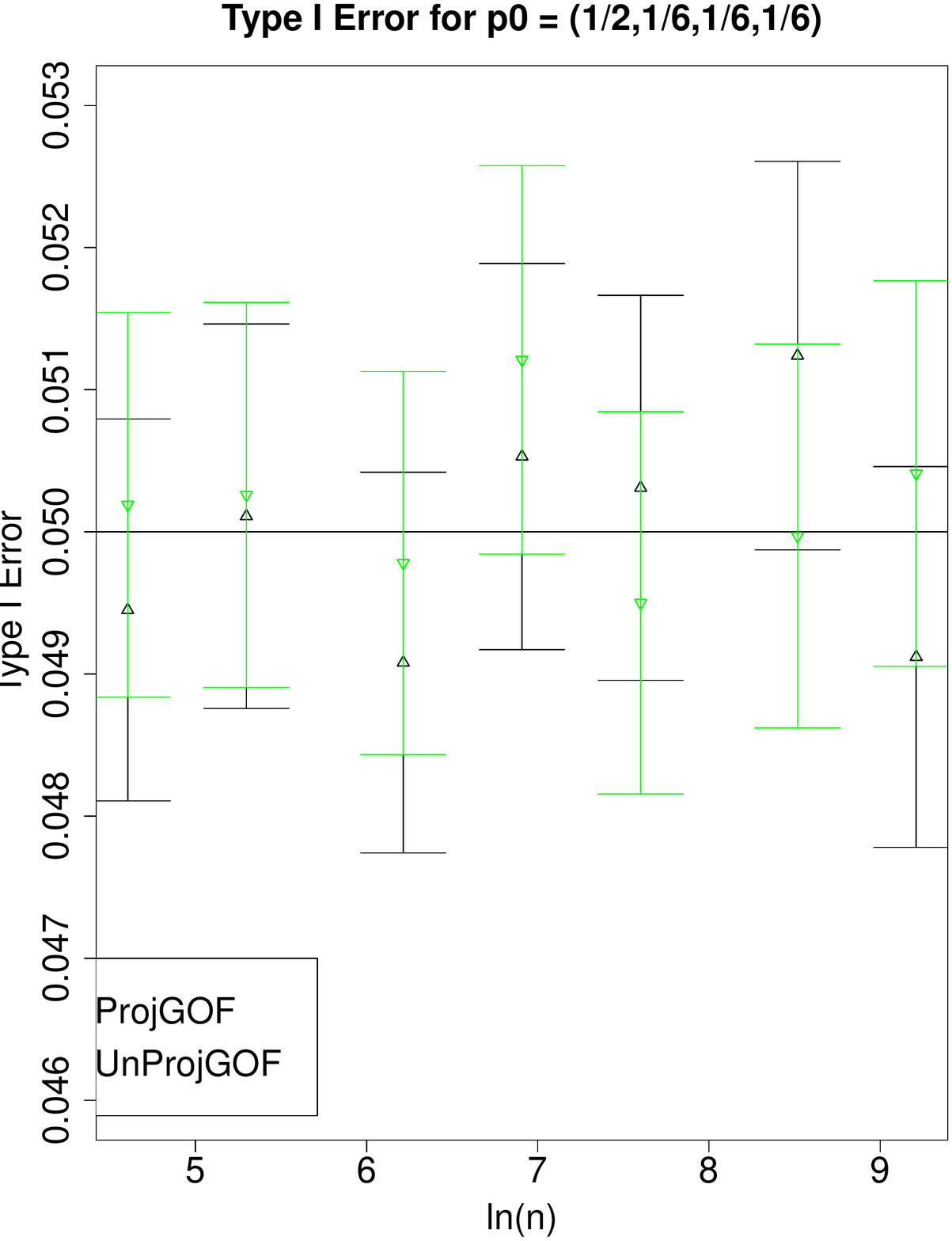}
\end{flushright}
\end{subfigure}
\caption{Empirical Type I Error for our new goodness of fit tests in $\GOF$ with error bars corresponding to 1.96 times the standard error in $100,000$ trials. We set $\rho = 0.001$ which corresponds to a variance of $1,000$ for the additional noise to the counts due to privacy.  The horizontal line corresponds to the target $\alpha = 0.05$ Type I error that we permit.  \label{fig:GOF_signif} }
\end{center}
\end{figure}
\fi

\ifnum\aistats=1
\begin{figure}
\centering
\includegraphics[scale=0.25]{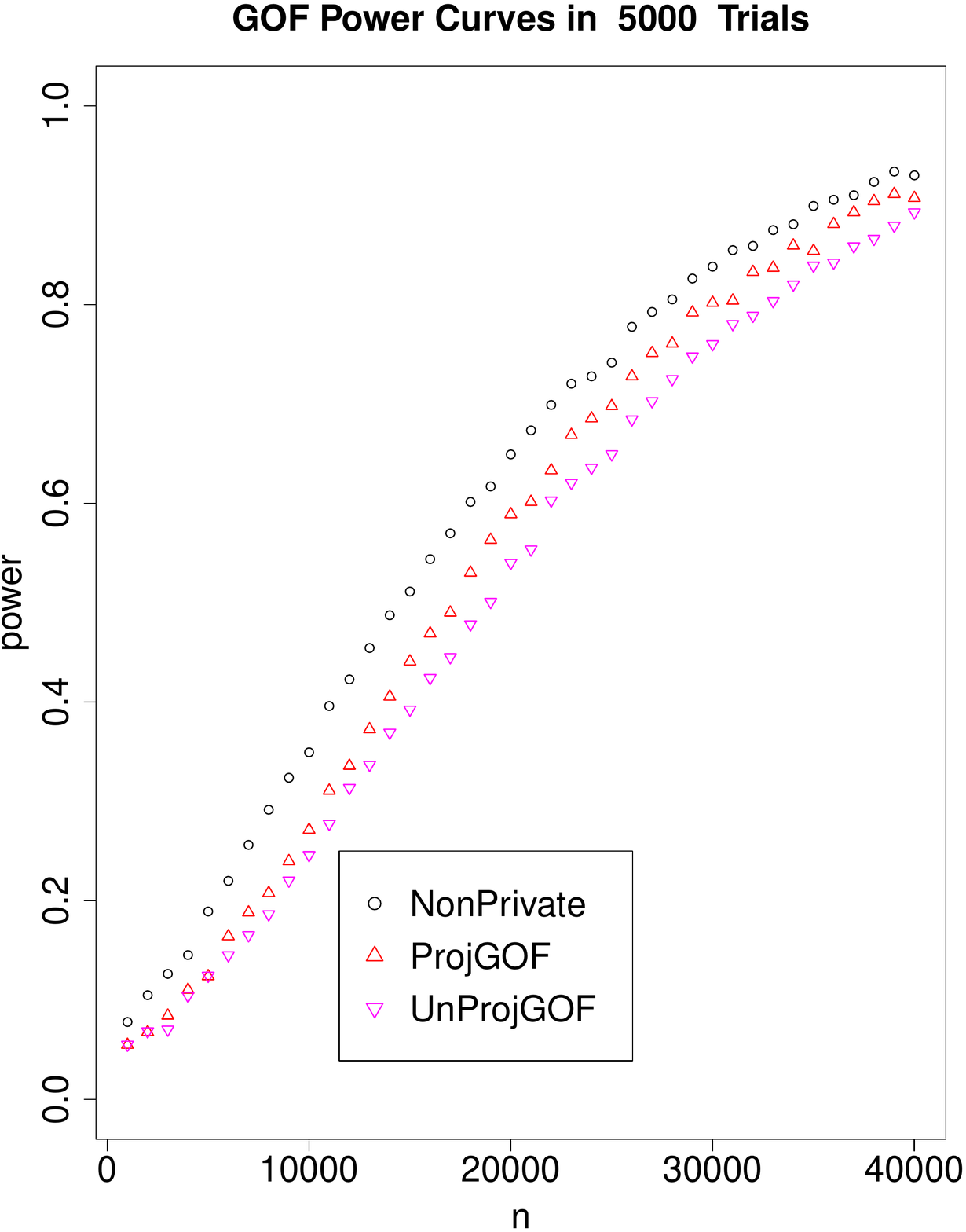}
\caption{\small Comparing power between the projected and unprojected statistics in $\GOF$ with the classical nonprivate test with $5000$ trials each, $\rho =0.001$ and $\alpha = 0.05$. \label{fig:power_gof}}
\end{figure}
\fi

We then empirically check the power of our new tests in $\GOF$ for both the projected and unprojected statistic.  Subject to the constraint that our tests achieve Type I error at most $\alpha$, we seek to maximize \emph{power}, or the probability of rejecting the null hypothesis when a distribution $\bp^1 \neq \nullbp$, called the \emph{alternate hypothesis}, is true.  
We expect to see the projected statistic achieve higher power than the unprojected statistic due to \Cref{thm:compare}.  Further, the fact that the critical value we use for the projected statistic is smaller than the critical value for the unprojected statistic might lead to the projected statistic having higher power.  

Here we present a typical experimental scenario. For our experiments, we set the null hypothesis $\nullbp = (1/2,1/6,1/6,1/6)$ and alternate hypothesis $\bp^1 = \nullbp + 0.01 \cdot (1,-1/3,-1/3,-1/3)$ for various sample sizes (we empirically found this to be a tough alternative hypothesis for our statistics).  For each sample size $n$, we sample $5,000$ independent datasets from the alternate hypothesis and test $H_0: \bp = \nullbp$ in $\GOF$.  We present the resulting power plots in \Cref{fig:power_gof} for $\GOF$ from \Cref{alg:gof}.
We label  ``NonPrivate" as the classical chi-square goodness of fit test used on the actual data (and thus not private).  Further, we write ``ProjGOF" as the test from $\GOF$ with the projected statistic whereas ``UnProjGOF" uses the unprojected statistic. Clearly in our results the projected outperforms the unprojected statistic.   


\ifnum\aistats=1
\begin{figure}
\centering
\includegraphics[scale=0.25]{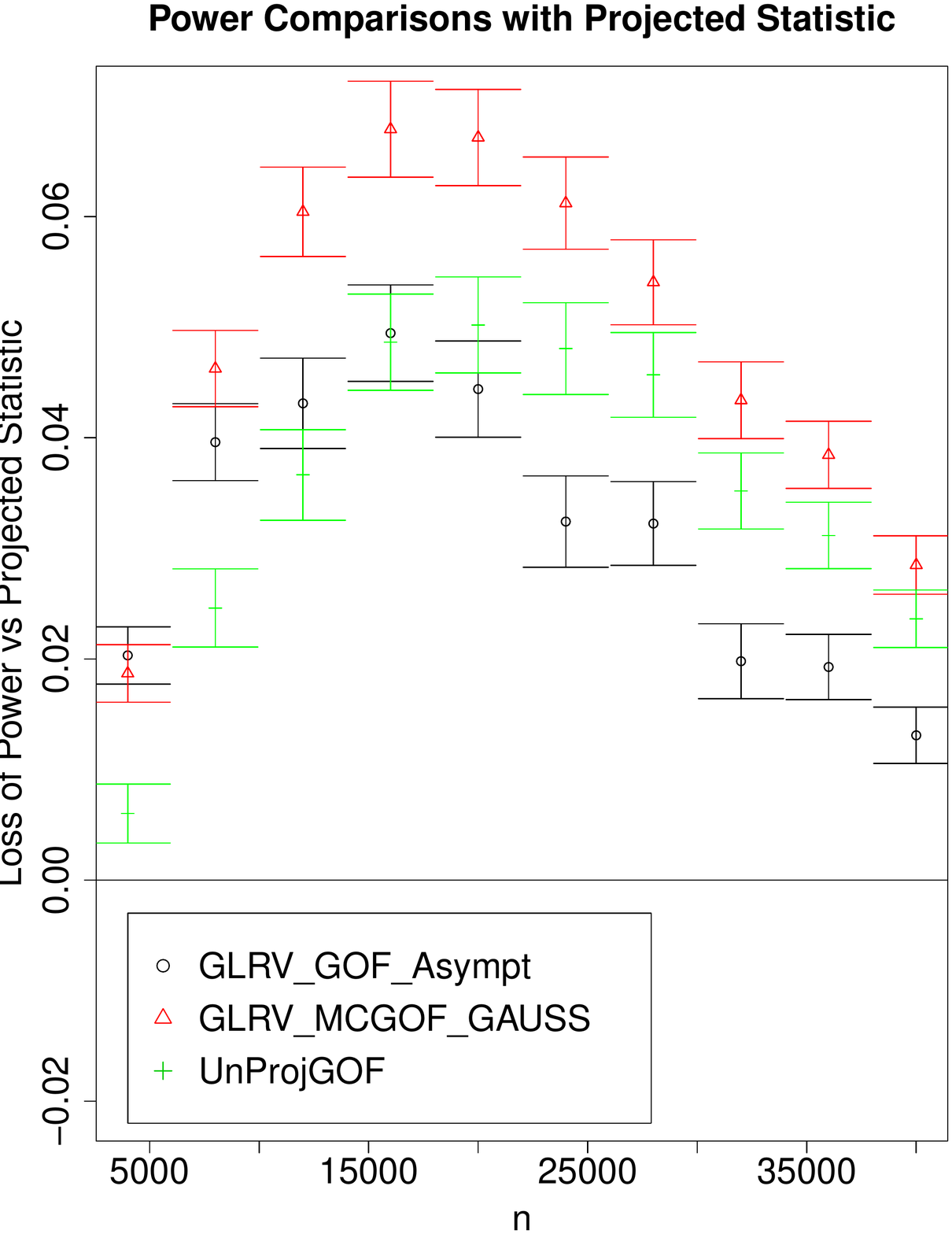}
\caption{\small The empirical power loss from using other private goodness of fit tests instead of the projected statistic in $\GOF$ 
for $100,000$ trials, $\rho =0.001$ and $\alpha = 0.05$. \label{fig:power_gof_compare}}
\end{figure}
\fi

We then compare the projected and unprojected statistic in $\GOF$ to prior work in \Cref{fig:power_gof_compare}. Since the projected statistic outperforms the other tests, we plot the difference in power between the projected statistic and the other tests.
We label ``GLRV\_MCGOF\_GAUSS" as the Monte-Carlo (MC) based test with Gaussian noise from \cite{GLRV16},\footnote{We set the the number of MC trials $m = 59$ in these experiments, which guarantees at most $5\%$ Type I error.} and ``GLRV\_GOF\_Asympt" as the hypothesis test based on the asymptotic distribution with Gaussian noise from \cite{GLRV16,WLK15}.  
Note that the error bars show $1.96$ times the standard error in the difference of proportions from $100,000$ trials, giving a $95\%$ confidence interval.

\ifnum\aistats=0
\begin{figure}
\begin{center}
\begin{subfigure}{.45\textwidth}
\includegraphics[width=\linewidth]{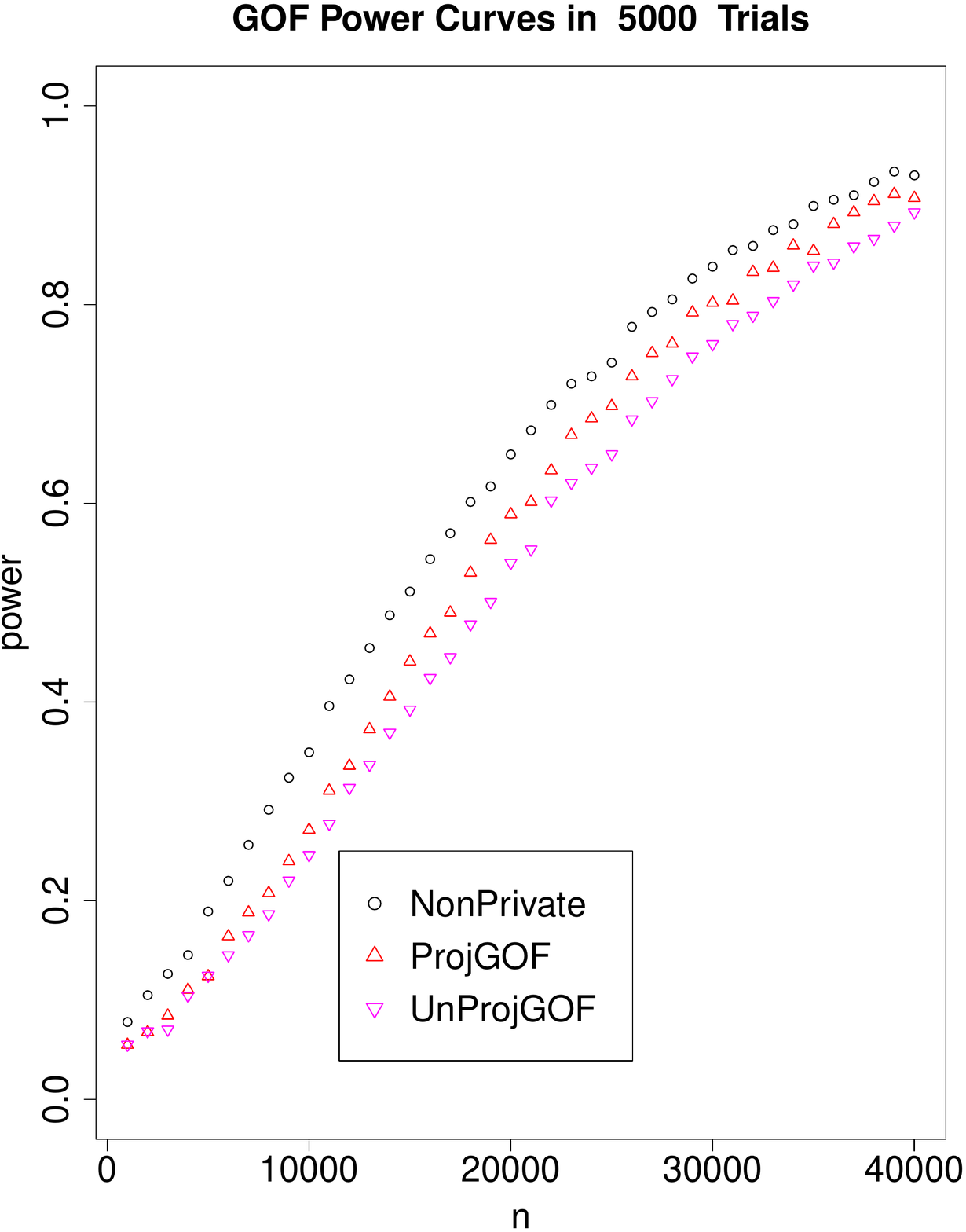}
\begin{flushleft}
\caption{A comparison of power for goodness of fit testing between the projected and unprojected statistics in $\GOF$ with the classical nonprivate test for various $n$ with $5,000$ trials each.  \label{fig:power_gof}}
\end{flushleft}
\end{subfigure}
\hspace{10mm}
\begin{subfigure}{.45\textwidth}
\begin{flushright}
\includegraphics[width=\linewidth]{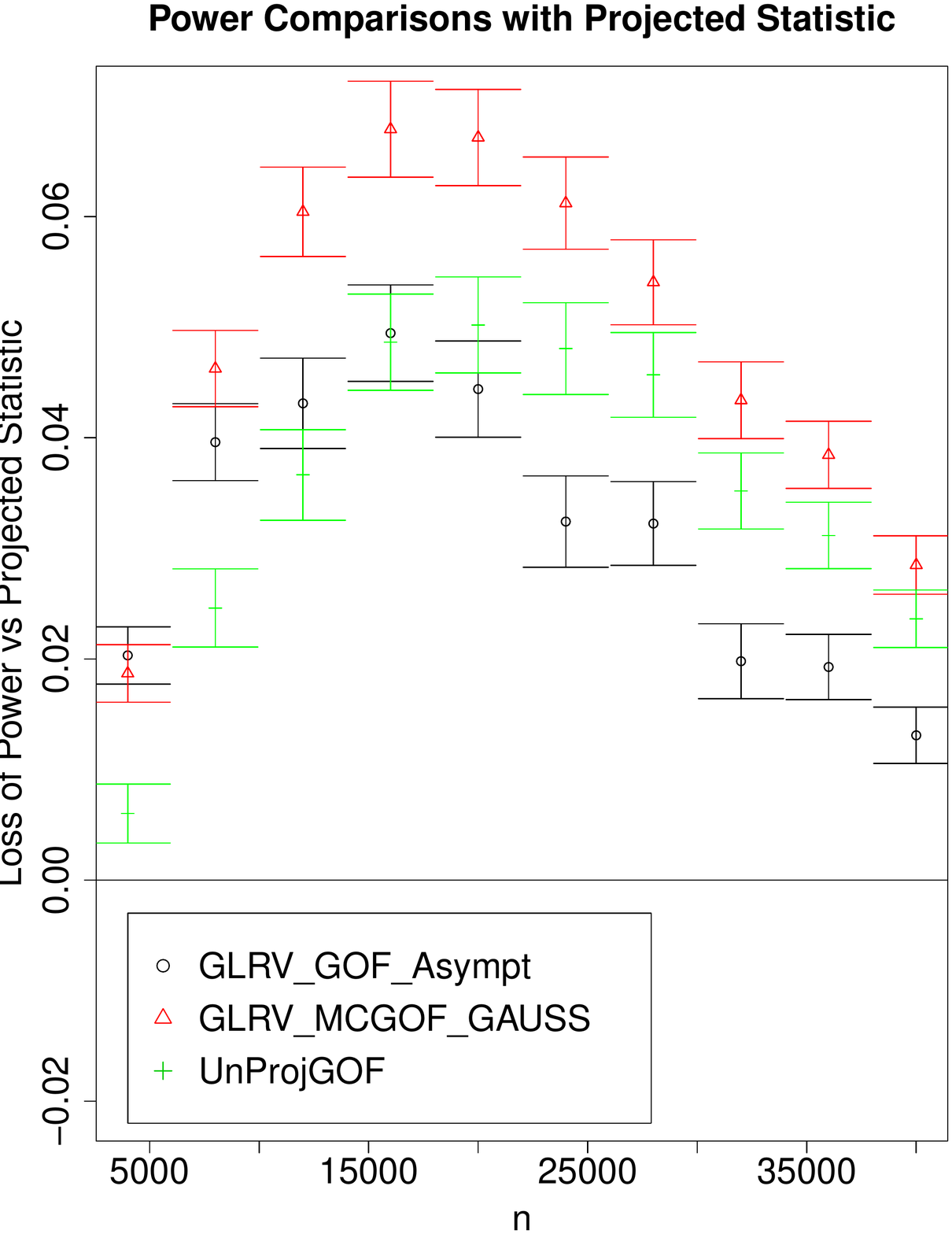}
\caption{The empirical loss in power from using other private goodness of fit tests instead of using the projected statistic in $\GOF$ with error bars corresponding to $1.96$ times the standard error of each difference for $100,000$ trials. \label{fig:power_gof_compare}}
\end{flushright}
\end{subfigure}
\caption{Empirical power results for our new goodness of fit tests in $\GOF$ with $\alpha = 0.05$ and a comparisons to previous private tests in \cite{GLRV16}.  We use $\rho = 0.001$, which corresponds to the variance of $1,000$ for the additional noise to the counts. }
\end{center}
\end{figure}
\fi


\section{General Chi-Square Private Tests}\label{sec:indep}
We now consider the case where the null hypothesis contains many distributions, 
so that the best fitting distribution must be estimated and used in the test statistics.
The data is multinomial $\XX{n} \sim \mult(n, \bp(\truetheta))$ and $\bp$ is a function that converts parameters into a $d$-dimensional multinomial probability vector. The null hypothesis is
$H_0: \truetheta \in \Theta$; i.e. $\bp(\truetheta)$ belongs to a subset of a lower-dimensional manifold.
We again use Gaussian noise $Z \sim N(0,1/\rho\cdot I_d)$ to ensure $\rho$-zCDP, and we define
\begin{equation}
\pUU{\rho}(\theta)  \stackrel{\defn}{=} \sqrt{n} \left( \frac{\XX{n}+ Z}{n} - \bp\left(\theta\right) \right).
\label{eq:covar_comp}
\end{equation}
With $\truetheta$ being the unknown true parameter,
we are now ready to define our two test statistics in terms of some function $\phi: \R^d \to \R$, such that $\phi(\XX{n} + Z) \stackrel{P}{\to} \truetheta$ (recall from Section \ref{subsec:minchi} that $\phi$ is a simple but possibly a suboptimal estimate of the true parameter $\truetheta$ based on the noisy data) and the covariance matrix 
\mathmath
\PrivCovar{\rho}(\theta) \stackrel{\defn}{=} \diag\left(\bp(\theta) \right) - \bp(\theta)\bp(\theta)^\intercal + 1/\rho \cdot I_d.
\mathmath

We define the \emph{unprojected} statistic $\RR{n}_{\rho}(\theta)$ as follows:
\begin{align}
\hat M \stackrel{\defn}{=} &\left(\PrivCovar{n \rho}\left(\phi(\XX{n} + Z)\right)\right)^{-1}\nonumber\\
\RR{n}_{\rho}(\theta) \stackrel{\defn}{=} &~\pUU{\rho}(\theta)^\intercal  \hat M  \pUU{\rho}(\theta).
\label{eq:R_unproj}
\end{align}
This is a specialization of \eqref{eq:R} in \Cref{subsec:minchi} with the following substitutions:
$\VV{n}= \left( \frac{\XX{n} + Z }{n}\right)$, $A(\theta) = \bp(\theta)$, and $M(\theta) = \left(\PrivCovar{n \rho}(\theta)\right)^{-1}$.  

For the \emph{projected} statistic $\hRR{n}_{\rho}(\theta)$,
the corresponding substitutions are
$\Proj=I_d - \frac{1}{d} \one\one^\intercal$, $\VV{n}= \Proj\cdot\left( \frac{\XX{n} + Z }{n}\right)$, $A(\theta) =\Proj \cdot \bp(\theta)$, and again $M(\theta) = \left(\PrivCovar{n\rho}(\theta)\right)^{-1}$ 
giving:
\begin{equation}
  \hRR{n}_{\rho}(\theta) \stackrel{\defn}{=}\pUU{\rho}(\theta)^\intercal  \cdot \Proj  \hat M \Proj \cdot  \pUU{\rho}(\theta).
\label{eq:R_proj}
\end{equation}
We then assume that for both the projected and unprojected statistic \Cref{assumpt:mean} holds using their relative vectors $\VV{n}$, $A(\theta)$, and matrix $M(\theta)$.  We now present the asymptotic distribution of both statistics, which is proved using the result in \Cref{thm:ferg} \ifnum\aistats=1 
 $~$(the full proof is in the supplementary file)\fi.

\begin{theorem}
Under $H_0:\truetheta \in \Theta$, the following are true as $n\rightarrow\infty$. Setting $\mintheta=\arg\min_{\theta\in \Theta} \RR{n}_{\rho_n}(\theta)$ we have 
$
\RR{n}_{\rho_n}(\mintheta) \stackrel{D}{\to} \chi^2_{ d - k } 
$
if $n\rho_n\to\rho>0$. Furthermore, setting $\mintheta=\arg\min_{\theta\in \Theta} \hRR{n}_{\rho_n}(\theta)$ we have 
$\hRR{n}_{\rho_n}(\mintheta) \stackrel{D}{\to} \chi^2_{d-k-1}$ if $n\rho_n\rightarrow\rho$ or $n\rho_n\rightarrow\infty$.
\label{thm:comp_main}
\end{theorem}
\ifnum\aistats=0
\begin{proof}
To prove this result, we appeal to \Cref{thm:ferg}.  For the unprojected statistic $\RR{n}_\rho(\cdot)$ we have that $C(\theta) = \Sigma_\rho(\theta)$ and the middle matrix $M(\theta)$ is simply the inverse of it, which satisfies the hypotheses of \Cref{thm:ferg}.  

For the projected statistic $\hRR{n}_{n\rho_n}(\cdot)$, we will write $C(\theta) = \Proj\cdot\PrivCovar{n\rho_n}(\theta)\cdot \Proj$, $M(\theta) =  \PrivCovar{n\rho_n}^{-1}(\theta)$, and $\dot{A}(\theta) = \Proj \cdot \nabla \bp(\theta) \in \R^{d \times k}$.   Note that $C(\theta)$ has rank $d-1$ for all $\theta \in \Theta$ in a neighborhood of $\truetheta$ and all $n$.  We will now show that we can satisfy the hypotheses in \Cref{thm:ferg} with these matrices, i.e. we show the following two equalities hold for all $\theta \in \Theta$
\begin{align*}
C(\theta) \cdot M(\theta) \cdot C(\theta) = C(\theta) \qquad \& \qquad C(\theta) \cdot M(\theta)\cdot  \dot{A}(\theta) = \dot{A}(\theta).
\end{align*}
We first focus on proving the first equality $C(\theta) \cdot M(\theta) \cdot C(\theta) = C(\theta)$.  From \eqref{eq:middle_rho}, we can simplify the left hand side of the equality significantly by rewriting it as
$$
 \Proj \PrivCovar{n\rho_n}(\theta) \Proj - \frac{n \rho_n}{d} \Proj \PrivCovar{n\rho_n}(\theta) \pmb{1}\pmb{1}^\intercal  \PrivCovar{n\rho_n}(\theta)  \Proj
$$
We now show that $\Proj \PrivCovar{n\rho_n}(\theta) \pmb{1}\pmb{1}^\intercal  \PrivCovar{n \rho_n}(\theta) = 0$ for all $n$, which would prove this equality.  Note that $\PrivCovar{n\rho_n}(\theta)$ is symmetric and has eigenvector $\pmb{1}$ with eigenvalue $\frac{1}{n\rho_n}$.  Thus, 
$$
\Proj \PrivCovar{\rho}(\theta) \pmb{1}\pmb{1}^\intercal  \PrivCovar{n\rho_n}(\theta)  =\frac{1}{n^2\rho_n^2} \Proj\pmb{1}\pmb{1}^\intercal  = 0 \quad \forall n .
$$

We now prove the second equality $C(\theta) \cdot M(\theta)\cdot  \dot{A}(\theta) = \dot{A}(\theta)$.  We again use \eqref{eq:middle_rho} to simplify the left hand side of the equality:
\begin{align*}
\Proj & \PrivCovar{n \rho_n}(\theta) \left[  \PrivCovar{n \rho_n}^{-1}(\theta) - \frac{n \rho_n}{d} \cdot \pmb{1}\pmb{1}^\intercal \right] \nabla \bp(\theta) \\
& = \Proj  \left[ I_d - \frac{n \rho_n}{d}\cdot \PrivCovar{n \rho_n}(\theta)\pmb{1}\pmb{1}^\intercal \right] \nabla \bp(\theta) \\
& =  \Proj  \Proj \nabla \bp(\theta) \\
& = \Proj \nabla \bp(\theta).
\end{align*}

This completes the proof for both cases $n \rho_n \to \rho>0$ and $n \rho_n \to \infty$.

\end{proof}
\fi
Again, the projected statistic has the same distribution under both private asymptotic regimes and matches the non-private chi-square test asymptotics.  We present our more general test $\MIN$ in \Cref{alg:min}.
The quick-and-dirty estimator $\phi(\cdot)$ is application-specific
(Section \ref{sect:IND} gives independence testing as an example).\footnote{For goodness-of-fit testing, $\phi$ always returns $\nullbp$ and $k=0$ so $\MIN$ is a generalization of $\GOF$.}  Further, for neighboring histogram data, we have the following privacy guarantee.

\begin{theorem}
 $\MIN(\cdot; \rho,\alpha,\phi,\Theta)$ is $\rho$-zCDP.
 \end{theorem}
%
\begin{algorithm}
\caption{ zCDP General Chi-Square Test}
\label{alg:min}
\begin{algorithmic}
\Procedure{\MIN}{$\XX{n}$; $\rho$, $\alpha$, $\phi$, $H_0: \truetheta \in \Theta$}
\\ Set $\pXX{n}\gets \XX{n} + Z$ where $Z \sim N(0,1/\rho \cdot I_d)$.
\\ Set $\hat M = \PrivCovar{n \rho}\left(\phi(\pXX{n}) \right)^{-1}$
\State \textbf{For the \emph{unprojected statistic}:}
\begin{align*}
 \stat(\theta) 
 = \frac{1}{n}\left(\pXX{n} - n\bp(\theta)  \right)^\intercal  \hat M \left(\pXX{n} - n\bp(\theta) \right)
\end{align*}
\State Set $\mintheta = \myargmin_{\theta \in \Theta} \stat(\theta)$
\State $t\gets (1-\alpha)$ quantile of $\chi^2_{d-k}$
\State \textbf{For the \emph{projected statistic}:}
\begin{align*}
 \stat(\theta) 
 = \frac{1}{n}\left(\pXX{n} - n\bp(\theta)\right)^\intercal  \Proj \hat M \Proj \left(\pXX{n} - n\bp(\theta)  \right)
\end{align*}
\State Set $\mintheta = \myargmin_{\theta \in \Theta} \stat(\theta)$
\State $t\gets (1-\alpha)$ quantile of $\chi^2_{d-k-1}$
\\ {\bf if } $\stat(\mintheta) > t$ {\bf then }
   Reject
\EndProcedure
\end{algorithmic}
\end{algorithm}

\subsection{Application - Independence Test \label{sect:IND}}
We showcase our general chi-square test $\MIN$ by giving results for independence testing. Conceptually, it is convenient to think of the data histogram as an $r\times c$ table, with $\bp_{i,j}$ being the probability a person is in the bucket in row $i$ and column $j$. We then consider two multinomial random variables $Y \sim \mult(1,\pipi{1})$ for $\pipi{1} \in \R^{r}$ (the marginal row probability vector) and $Y' \sim \mult(1,\pi^{(2)})$ for $\pipi{2} \in \R^{c}$  (the marginal column probability vector). Under the null hypothesis of independence between $Y$ and $Y^\prime$, $\bp_{i,j}=\pi^{(1)}_i\pi^{(2)}_j$. 
Generally, we write the probabilities as $\bp(\pi^{(1)}, \pi^{(2)}) = \pipi{1} \left(\pipi{2}\right)^\intercal $ so that 
\mathmath
X^{(n)} \sim \mult\left(n,\bp(\pi^{(1)}, \pi^{(2)})\right).
\mathmath
Thus we have the underlying parameter vector $\truetheta = \left(\pipi{1}_{1},\cdots,\pipi{1}_{r-1}, \pipi{2}_{1}, \cdots,\pipi{2}_{c-1} \right)$ - we do not need the last component of $\pipi{1}$ or $\pipi{2}$ because we know that each must sum to 1.  Also, we have $d = r c$ and $k = (r-1)+(c-1)$ in this case.  We want to test whether $Y$ is independent of $Y'$.  For our data, we are given a collection of $n$ independent trials of $Y$ and $Y'$.  We then count the number of joint outcomes in a contingency table given in \Cref{table:contingency}.  Each cell in the contingency table contains element $\XX{n}_{i,j}$ that gives the number of occurrences of $Y_i = 1$ and $Y'_j =1$. 
Since our test statistics notationally treat the data as a vector, when needed, we convert $\XX{n}$ 
to a vector that goes from left to right along each row of the table. 

\begin{table}[ht]
\caption{Contingency Table.} 
\centering 
\begin{tabular}{ c | c | c |c | c | c}
\hline
$Y \quad \backslash \quad Y'$ & 1 & 2 & $\cdots$ & $c$ & Marginals \\
\hline
1 & $\XX{n}_{1,1}$ & $\XX{n}_{1,2}$ & $\cdots$ & $\XX{n}_{1,c}$ & $\XX{n}_{1,\cdot}$ \\
\hline
2 &$\XX{n}_{2,1}$ &$\XX{n}_{2,2}$ & $\cdots$ & $\XX{n}_{2,c}$ & $\XX{n}_{2,\cdot}$ \\
\hline
$\vdots$ &$\vdots$ & $\vdots$& $\ddots$& $\vdots$ & $\vdots$ \\
\hline
$r$ &$\XX{n}_{r,1}$ &$\XX{n}_{r,2}$ &$\cdots$ &$\XX{n}_{r,c}$ & $\XX{n}_{r,\cdot}$ \\
\hline
Marginals & $\XX{n}_{\cdot, 1}$&$\XX{n}_{\cdot, 2}$ & $\cdots$ & $\XX{n}_{\cdot, c}$& $n$
\end{tabular}
\label{table:contingency}
\end{table}

In order to compute the statistic $\RR{n}_\rho(\mintheta)$ or $\hRR{n}_\rho(\mintheta)$ in $\MIN$, we need to find a quick-and-dirty estimator $\phi(\XX{n} + Z)$ that converges in probability to $\bp\left(\pipi{1},\pipi{2} \right)$ as $n \to \infty$.  We will use the estimator for the unknown probability vector based on the marginals of the table with noisy counts, so that for  na\"{i}ve estimates $\tilde \pi_i^{(1)} = \frac{\XX{n}_{i,\cdot} + Z_{i,\cdot}}{\tilde n}$, $\tilde \pi_j^{(2)} = \frac{\XX{n}_{\cdot, j}+Z_{\cdot,j}}{\tilde n}$ where $\tilde n = n + \sum_{i,j} Z_{i,j}$ we have\footnote{We note that in the case of small sample sizes, we follow a common rule of thumb where if any of the expected cell counts are less than $5$, i.e. if $n \tilde \pi_i^{(1)}\tilde \pi_j^{(2)}\leq 5$ for any $(i,j) \in [r]\times [c]$, then we do not make any conclusion.}

 \begin{equation}
 \phi\left(X^{(n)} + Z \right) =\left( \tilde\pi^{(1)}_1, \cdots,  \tilde\pi^{(1)}_{r-1},  \tilde\pi^{(2)}_1,\cdots,  \tilde\pi^{(2)}_{c-1}  \right)
 \label{eq:p_MLE}
 \end{equation}
\ifnum\aistats=1 
Note that $Z \sim N(0,1/\rho_n \cdot  I_{rc})$ so it is easy to see that under both private asymptotic regimes ($n\rho_n\rightarrow\rho$ and $n\rho_n\rightarrow\infty$)  we have $ \tilde\pi^{(1)}_i \stackrel{P}{\to} \pipi{1}_i$ and $\tilde\pi^{(2)}_j \stackrel{P}{\to} \pipi{j}_j$ for all $i \in [r]$ and $j \in [c]$ as $n\rightarrow\infty$.  
\else
Note that as $n \to \infty$, the marginals converge in probability to the true probabilities even for $Z \sim N(0,1/\rho_n \cdot  I_{rc})$ with $\rho_n = \omega(1/n^2)$, i.e. we have that $ \tilde\pi^{(1)}_i \stackrel{P}{\to} \pipi{1}_i$ and $\tilde\pi^{(2)}_j \stackrel{P}{\to} \pipi{j}_j$ for all $i \in [r]$ and $j \in [c]$.  Recall that in \Cref{thm:comp_main}, in order to guarantee the correct asymptotic distribution we require the $n \rho_n \to \rho>0$, or in the case of the projected statistic, we need $\rho_n  = \Omega(1/n)$.  Thus, \Cref{thm:comp_main} imposes more restrictive settings of $\rho_n$ for the unprojected statistic than what we need in order for the na\"ive estimate to converge to the true underlying probability.  For the projected statistic, we only need $\rho_n = \omega(1/n)$ to satisfy the conditions in \Cref{thm:comp_main} and for $\phi(\XX{n} + Z) \stackrel{P}{\to}\bp\left(\pipi{1},\pipi{2} \right).$ 
\fi

We then use this statistic $\phi(\XX{n}+Z)$ in our unprojected and projected statistic in $\MIN$ to have a $\rho$-zCDP hypothesis test for independence between two categorical variables.  Note that in this setting, the projected statistic has a $\chi^2_{(r-1)(c-1)}$ distribution, which is exactly the same asymptotic distribution used in the classical Pearson chi-square independence test.  

For our results we will again fix $\alpha = 0.05$ and $\rho = 0.001$. 
\ifnum\aistats=1
We verify experimentally in the supplementary file that our tests achieve at most $\alpha$ Type I error.
\else
In \Cref{fig:IND_signif} we give the empirical Type I error for our independence tests given in $\MIN$ for both the projected and unprojected statistic in $100,000$ trials for various $n$ and data distributions.  We note that for small sample sizes we are achieving much smaller Type I Errors than the target $\alpha$ due to the fact that sometimes the noise forces us to have small expected counts ($<5$ in any cell) in the contingency table based on the noisy counts, in which case our tests are inconclusive and fail to reject $H_0$.  
\fi

\ifnum\aistats=0
\begin{figure}
\begin{center}
\begin{subfigure}{.45\textwidth}
\begin{flushleft}
\includegraphics[width=\linewidth]{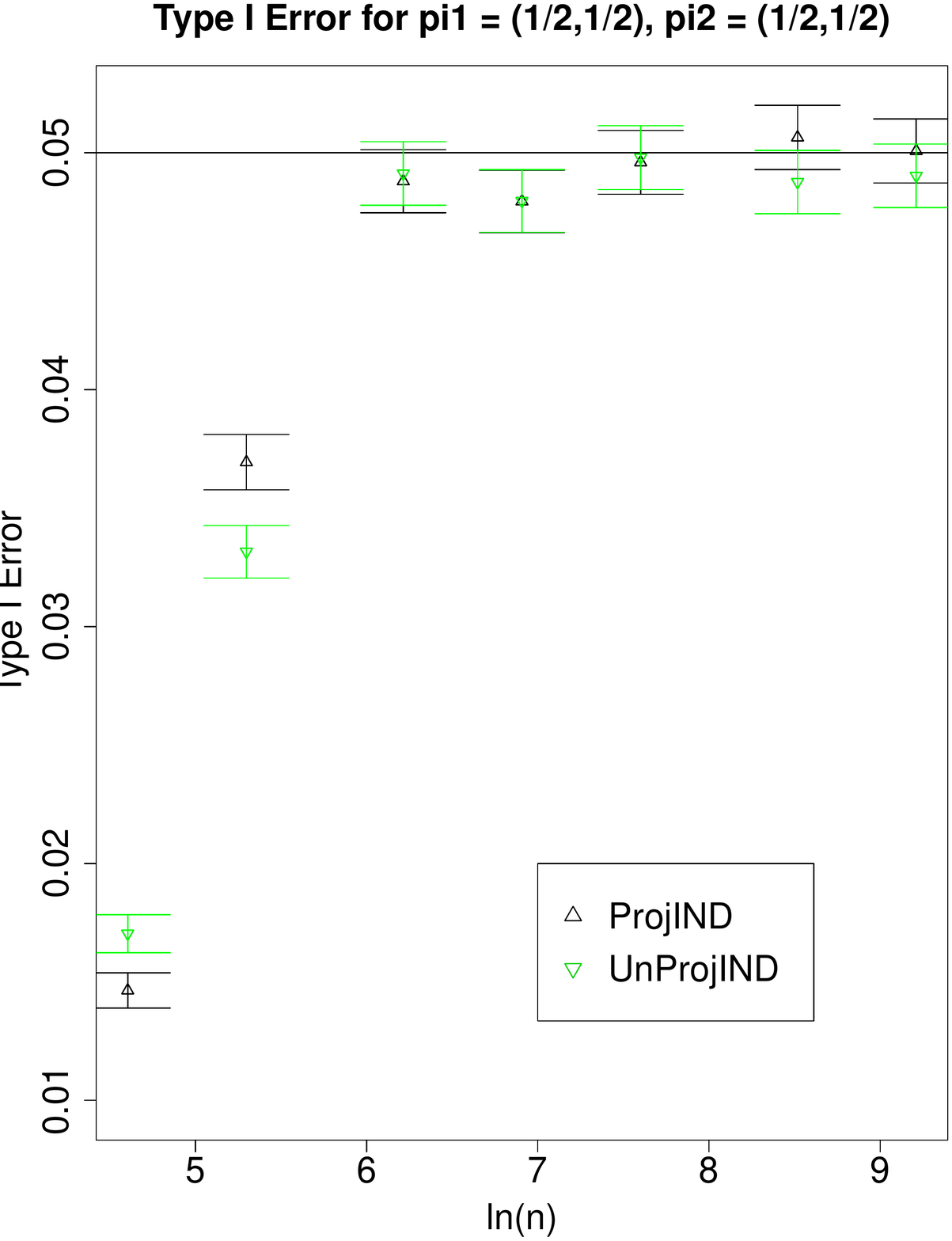}
\end{flushleft}
\end{subfigure}
\begin{subfigure}{.1\textwidth}
\end{subfigure}
\hspace{10mm}
\begin{subfigure}{.45\textwidth}
\begin{flushright}
\includegraphics[width=\linewidth]{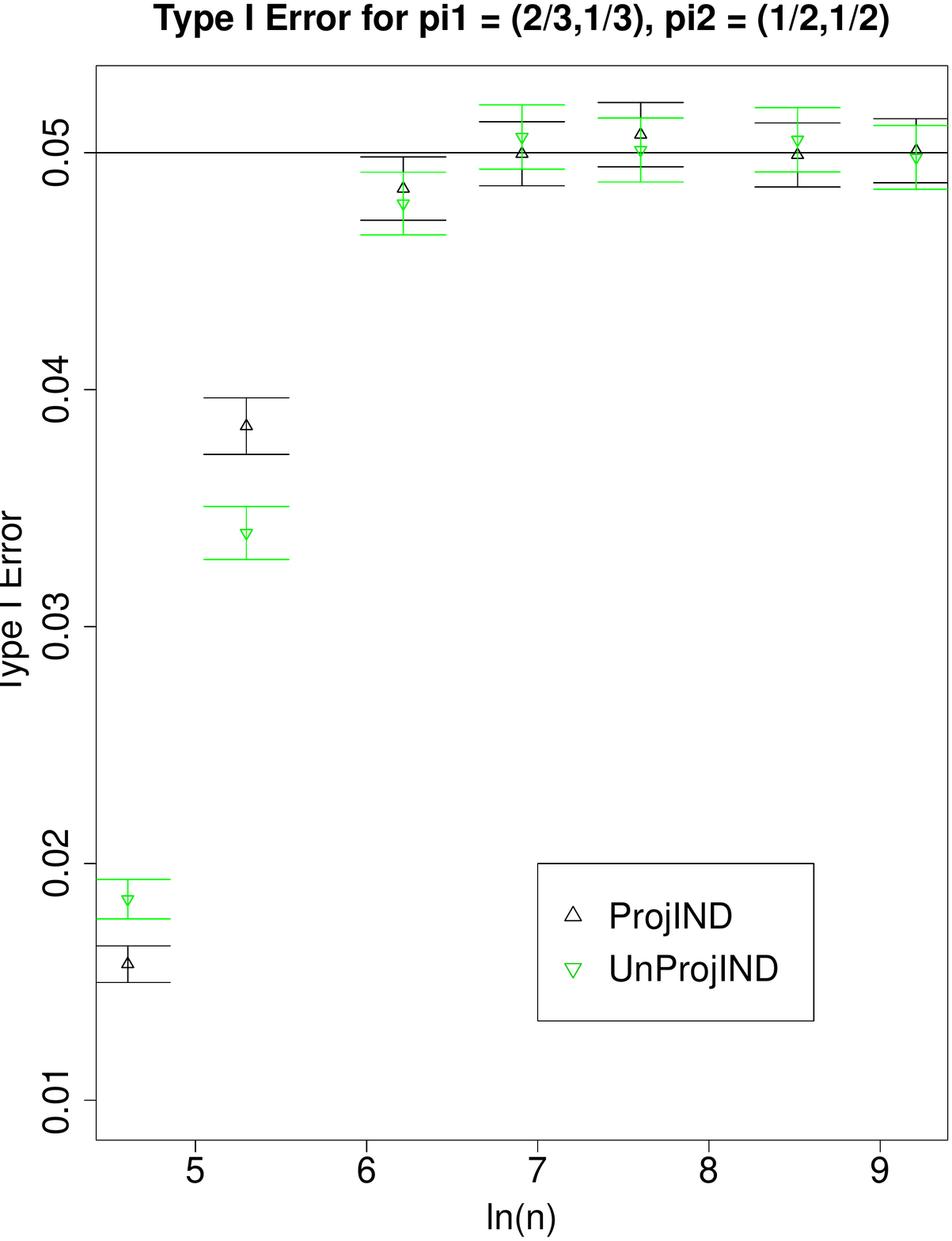}
\end{flushright}
\end{subfigure}
\caption{Empirical Type I Error for our new independence tests in $\MIN$ with 1.96 times the standard error in $100,000$ trials.  We set $\rho = 0.001$ which corresponds to variance $1,000$ due to noise in each cell count.  It is desired to have Type I error at most $\alpha=0.05$, which is given as the horizontal line.  \label{fig:IND_signif} }
\end{center}
\end{figure}
\fi


\ifnum\aistats=1
\begin{figure}
\begin{center}
\includegraphics[scale=0.25]{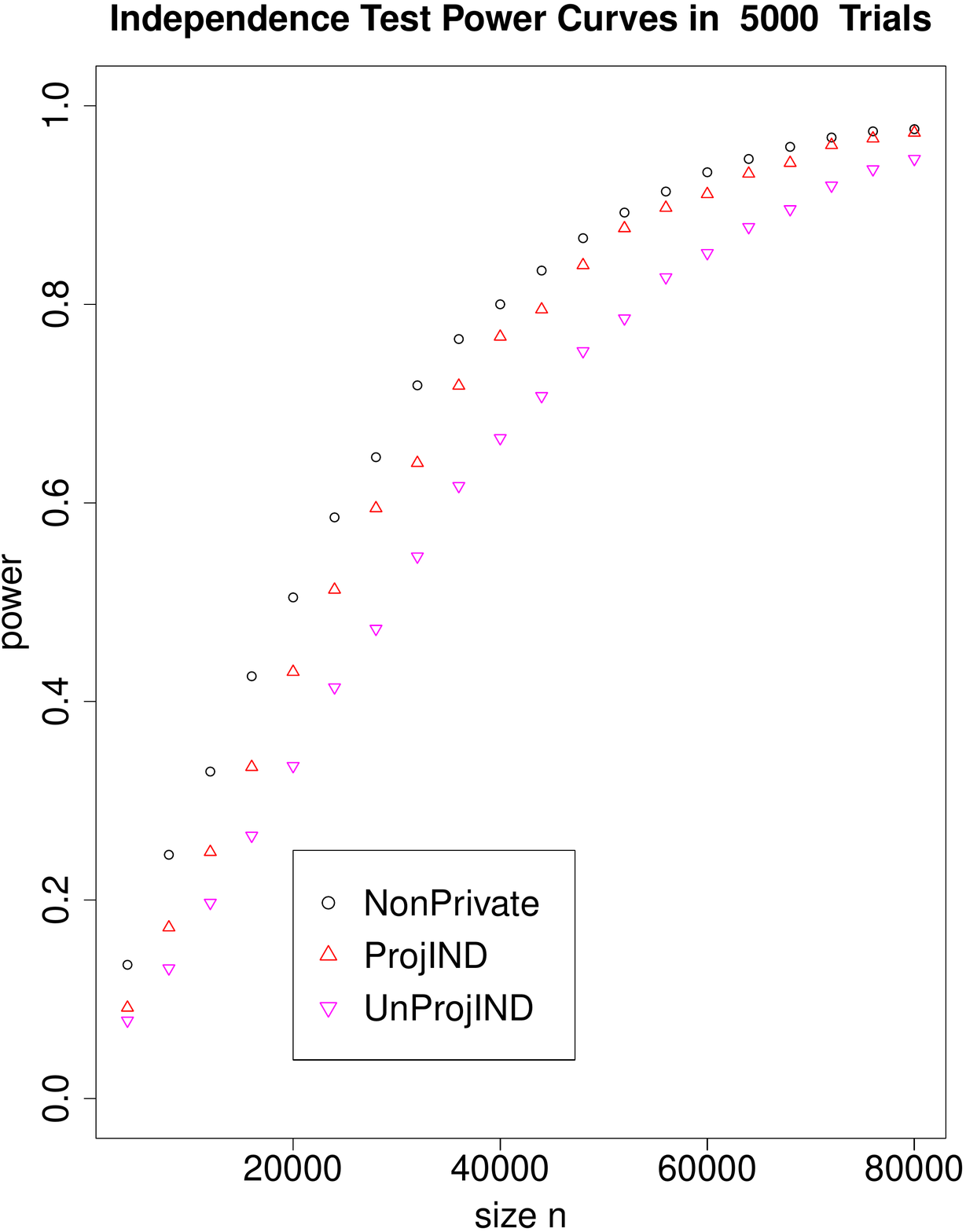}
\caption{\small Comparing power between the projected and unprojected statistics in $\MIN$ for independence testing with the classical Pearson chi-square test with $5000$ trials each, $\rho =0.001$ and $\alpha = 0.05$. \label{fig:ind_power}}
\end{center}
\end{figure}
\fi

We then compare the power $\MIN$ achieves for either of our test statistics.  As a sample of our experiments, we set $r = c = 2$ and $\pipi{1} =(2/3,1/3), \pipi{2} = (1/2,1/2)$.  We then sample our contingency table $X^{(n)}$ from $\mult(n, \bp(\pipi{1},\pipi{2}) + \pmb{\Delta})$ where $\pmb{\Delta} = 0.01 \cdot (1,0,-1,0)$, so that the null hypothesis is indeed false and should be rejected.  We give the empirical power of $\MIN$ in \Cref{fig:ind_power} using both the unprojected $\RR{n}_\rho(\mintheta)$ from \eqref{eq:R_unproj} and projected statistic $\hRR{n}_\rho(\mintheta)$ from \eqref{eq:R_proj} for $5,000$ independent trials and various sample sizes $n$.  Note that again we pick $\mintheta$ from \Cref{thm:min_chi} relative to the statistic we use.   We label ``NonPrivate" as the classical Pearson chi-square test used on the actual data and``ProjIND" as the test from $\MIN$ with the projected statistic whereas ``UnProjIND" uses the unprojected statistic.

\ifnum\aistats=1
\begin{figure}
\begin{center}
\includegraphics[scale=0.25]{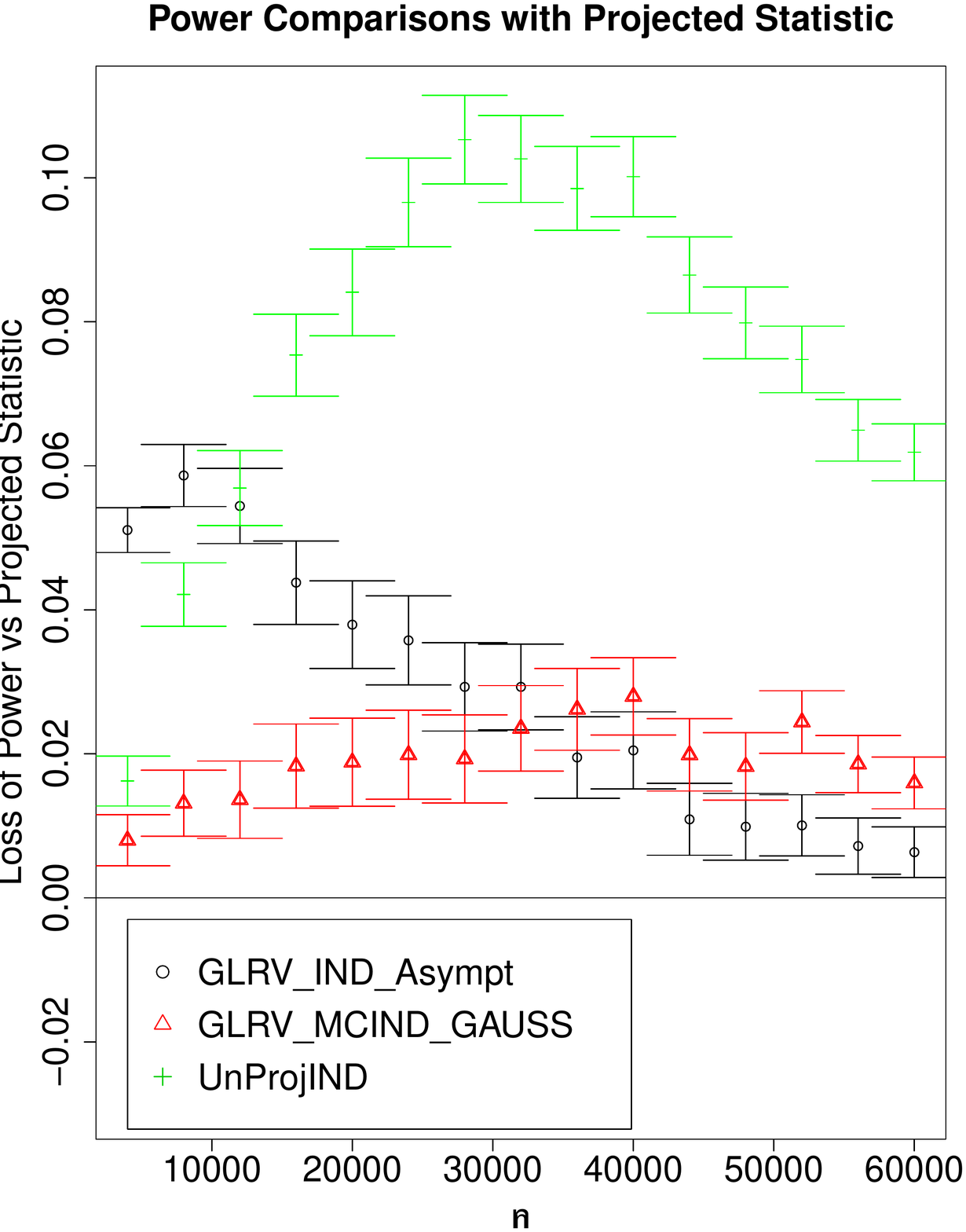}
\caption{\small The empirical power loss from using other private independence tests instead of the projected statistic in $\MIN$ 
for $50,000$ trials, $\rho =0.001$ and $\alpha = 0.05$. \label{fig:power_ind_compare}}
\end{center}
\end{figure}
\fi 

The projected statistic again outperforms prior work, so in \Cref{fig:power_ind_compare}, we plot the difference in power between the projected statistic in $\MIN$ and the competitors (the unprojected statistic and  independence tests from \cite{GLRV16}) in $50,000$ trials.  Note that we label ``GLRV\_MCIND\_GAUSS" as the MC based independence test with Gaussian noise and ``GLRV\_IND\_Asympt" as the hypothesis test based on the asymptotic distribution from \cite{GLRV16}.
\ifnum\aistats=1
Additional experiments can be found in the supplementary material.
\fi

\ifnum\aistats=0
\begin{figure}
\begin{center}
\begin{subfigure}{.45\textwidth}
\begin{flushleft}
\includegraphics[width=\linewidth]{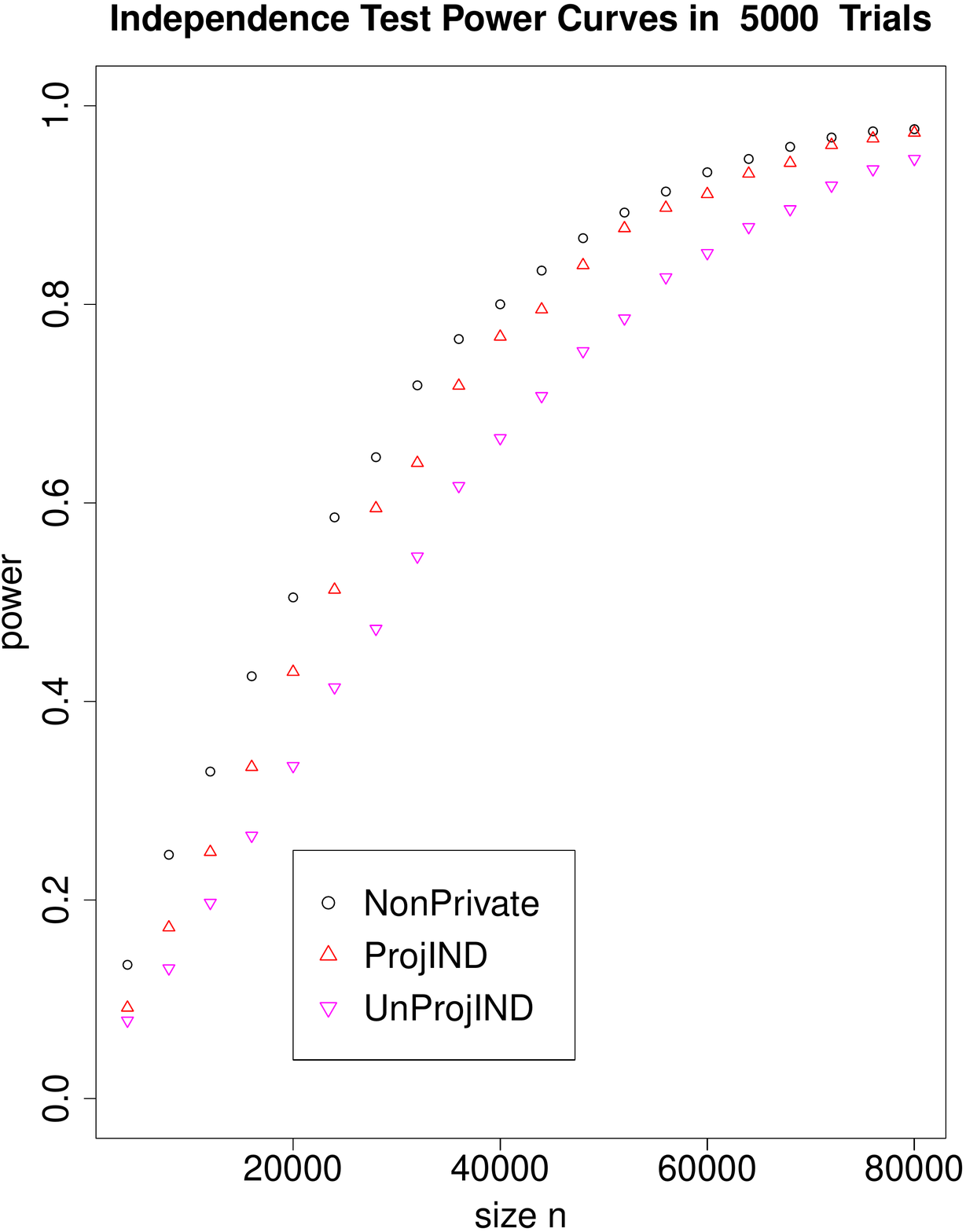}
\caption{A comparison of power for independence testing between the projected and unprojected statistics in $\MIN$ with the classical nonprivate test. We set $\alpha = 0.05$, $\rho = 0.001$ and $5,000$ to be the number of trials. \label{fig:ind_power}}
\end{flushleft}
\end{subfigure}
\hspace{10mm}
\begin{subfigure}{.45\textwidth}
\begin{flushright}
\includegraphics[width=\linewidth]{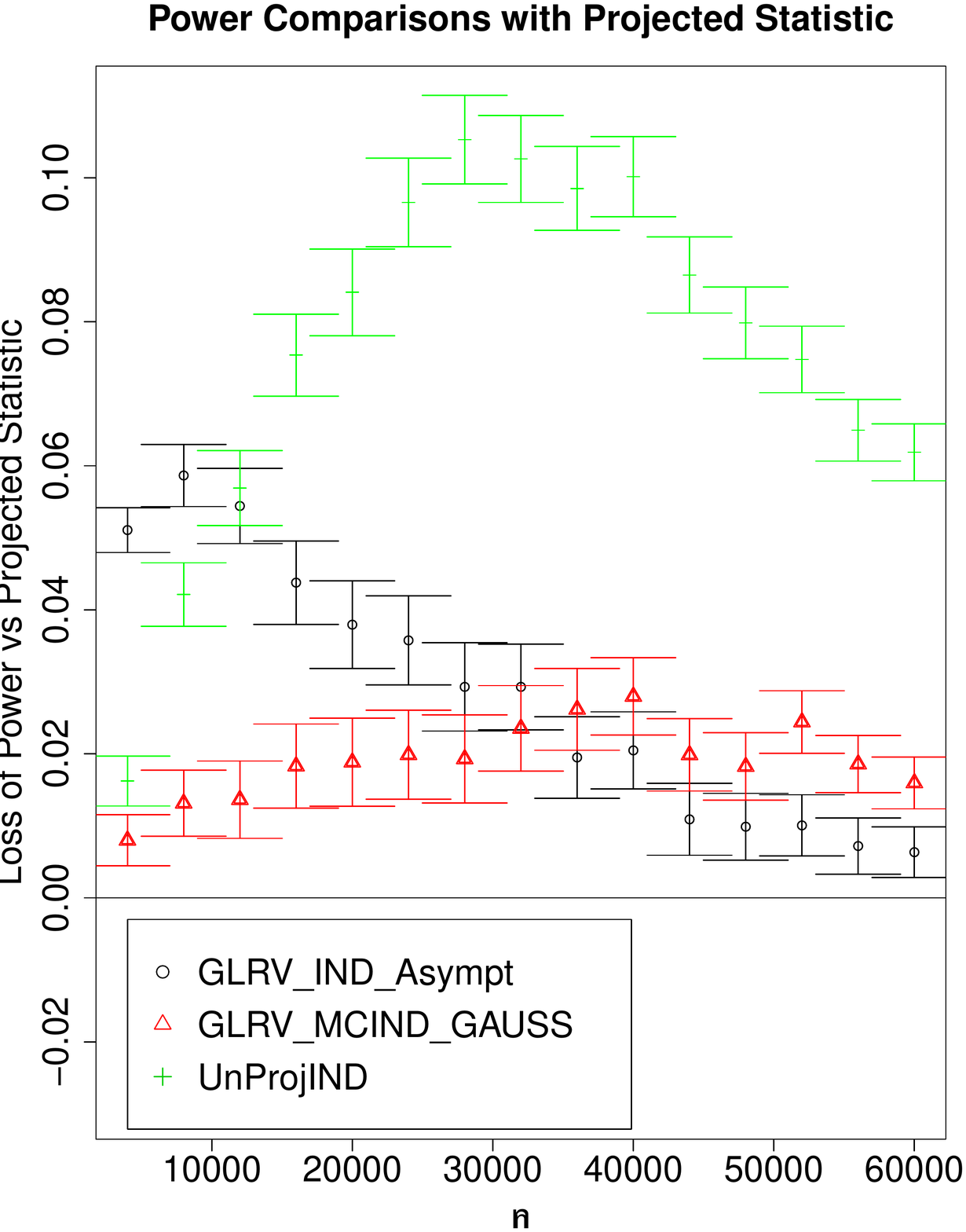}
\caption{The empirical power loss from using other private independence tests instead of using the projected statistic in $\MIN$ with error bars corresponding to $1.96$ times the standard error of each difference in $50,000$ trials. We set $\alpha = 0.05$ and $\rho= 0.001$.  \label{fig:power_ind_compare}}
\end{flushright}
\end{subfigure}
\caption{Empirical power results for our new independence tests in $\MIN$ and a comparisons to previous private tests in \cite{GLRV16} }
\end{center}
\end{figure}
\fi


\ifnum\aistats=0
\subsection{Application - GWAS Testing}
We next turn to demonstrating that our new general class of private hypothesis tests for categorical data significantly improves on existing private hypothesis tests even when extra structure is assumed about the dataset.  Specifically, we will be interested in GWAS data, which was the primary reason for why hypothesis tests for independence should be made private \cite{Homer08}.  We will then assume that $r = 3$ and $c=2$ and the data is evenly split between the two columns - as is the case in a control and case group.  For such tables, we can directly compute the sensitivity of the classical chi-square statistic $q(\cdot)$:
$$
q(\XX{n}) = \sum_{i=1}^3 \sum_{j = 1}^2 \frac{ n \cdot \left(\XX{n}_{i,j}- \frac{\XX{n}_{\cdot, j} \cdot \XX{n}_{i ,\cdot}}{n} \right)^2}{\XX{n}_{\cdot, j} \cdot \XX{n}_{i ,\cdot} } 
$$
\begin{lemma}[\cite{USF13,YFSU14}]
The $\ell_1$ and $\ell_2$ global sensitivity of the chi-square statistic $q(\cdot)$ based on a $3\times 2$ contingency table with positive margins and $n/2$
cases and $n/2$ controls is $\Delta(q) = 4n/(n+2)$.
\end{lemma}

Hence, a different approach for a private independence test is to add Gaussian noise with variance $\sigma^2  = \frac{\Delta_2(q)^2}{2\rho}$ to the statistic $q(\cdot)$ itself, which we call \emph{output perturbation}.  Our statistic is then simply Gaussian mechanism $\MG$ for statistic $q$.  We then compare the private statistic value with the distribution of $\cT_{\text{Gauss}}(n,\rho) = \chi^2_{2} + N\left(0,\sigma^2 \right)$ where the degrees of freedom is 2 because we have $(r-1)\cdot(c-1) = 2$.  Thus, given a Type I error of at most $\alpha$, we then set our critical value as $\tau_{\text{Gauss}}(\alpha;n,\rho)$ where
$$
 \prob{\cT_{\text{Gauss}}(n,\rho) > \tau_{\text{Gauss}}(\alpha;n,\rho)} = \alpha
$$
Hence, if $\MG(\XX{n})$ for the statistic $q$ is larger than $\tau_{\text{Gauss}}(\alpha;n,\rho)$ then we reject the null hypothesis.  

For our experiments, we again set $\rho = 0.001$ and $\alpha = 0.05$.  We fix the probability vector $(1/3,1/3,1/3)$ over the 3 rows in the first column whereas in the second column we set $(1/2,1/4,1/4)$, therefore the case and control groups do not produce the same outcomes.  In \Cref{fig:GWAS}, we show a comparison in the power between our test with the projected statistic, which assumes no structure on the data, and the output perturbation test, which crucially relies on the fact that the data is evenly split between the case and control groups.  We label ``ProjIND" and ``UnProjIND" as the tests from $\MIN$ with the projected statistic and unprojected statistic, respectively.  Further, we label ``YFSU\_Gauss" as the output perturbation tests for Gaussian noise proposed in \cite{YFSU14}.   Note that our new proposed test does significantly better than the output perturbation test, sometimes requiring $5$ times more samples to achieve the same level of power than for our projected statistic test.
\begin{figure}
\begin{center}
\includegraphics[width=.4\linewidth]{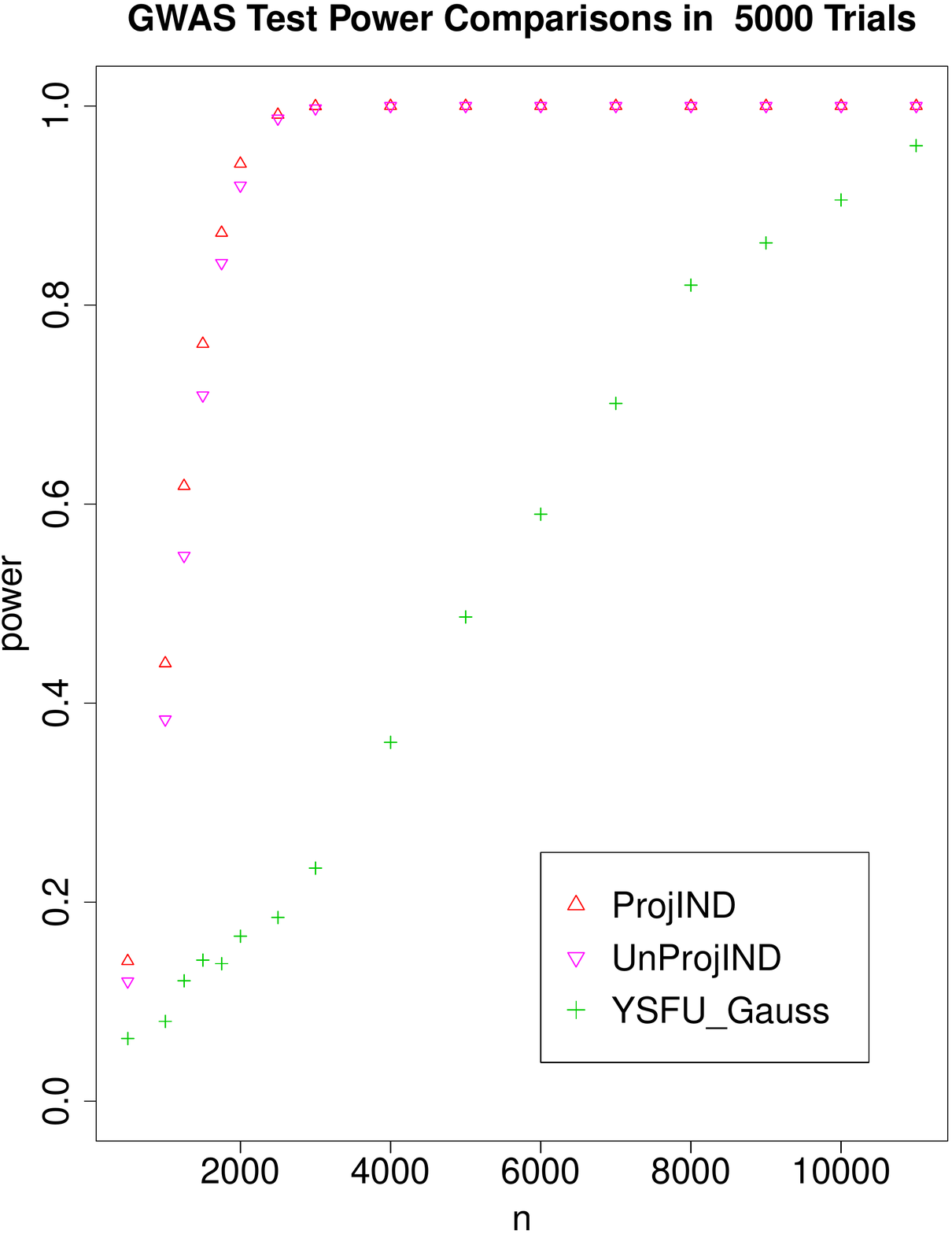}
\end{center}
\caption{A comparison of power between different hypothesis tests for independence testing for GWAS type datasets where the data is publicly known to be evenly split between the two columns and there are three rows in the contingency table.\label{fig:GWAS}}
\end{figure}
\fi

\ifnum\aistats=0
\section{General Chi-Square Tests with Arbitrary Noise Distributions}\label{sec:lapind}

We next show that we can apply our testing framework in \Cref{alg:min} for any type of noise distribution we want to include for privacy concerns.  For example, we consider adding Laplace noise rather than Gaussian noise if our privacy benchmark were (pure) differential privacy (DP).  In this case, we add Laplace noise with variance $8/\epsilon^2$ when computing the two statistics $\RR{n}_{\epsilon^2/8}(\mintheta)$ from \eqref{eq:R_unproj} and $\hRR{n}_{\epsilon^2/8}(\mintheta)$ from \eqref{eq:R_proj} so that the resulting tests will be $\epsilon$-DP and hence $\frac{\epsilon^2}{2}$-zCDP from \Cref{thm:reduction}.  Note that the resulting asymptotic distribution in this case will not be chi-square when we use noise other than Gaussian.  We will then rely on Monte Carlo (MC) sampling to find the critical value in which to reject the null hypothesis.  We give the MC based test which adds independent Laplace noise with variance $8/\epsilon^2$ in \Cref{alg:DP_MIN} and is thus $\epsilon$-DP, but any noise distribution can be used where we replace the parameter $1/\rho$ in the two statistics to be the variance of the noise that is added to each count.  In fact, Gaussian noise can be used in this framework although the asymptotic distribution seems to do well in practice even for small sample sizes.

\begin{algorithm}
\caption{DP Minimum Chi-Square Test using MC}
\label{alg:DP_MIN}
\begin{algorithmic}
\Procedure{$\texttt{DP-MC-MIN}$}{Histogram data $\XX{n} = (\XX{n}_1,\cdots, \XX{n}_d)$; $\epsilon$, $\alpha$, $H_0: \truetheta \in \Theta$, $m$ trials}
\State Set $\pXX{n}\gets \XX{n} + Z$ where $Z = (Z_1,\cdots, Z_d)$, where $Z_i \sim \text{Lap}(2/\epsilon)$.
\State Set $\hat M = \PrivCovar{n \rho}\left(\phi(\pXX{n}) \right)^{-1}$
\State \textbf{For the \emph{unprojected statistic}:}
\begin{align*}
 \stat(\theta) 
 = \frac{1}{n}\left(\pXX{n} - n\bp(\theta)  \right)^\intercal  \hat M \left(\pXX{n} - n\bp(\theta) \right)
\end{align*}
\State Set $\mintheta = \myargmin_{\theta \in \Theta} \stat(\theta)$
\State Sample $\{r_1,\cdots, r_m\}$ as  $m$ samples from the distribution of $\stat(\mintheta)$.
\State Set $\tau(\alpha,\epsilon)$ to be the $\lceil (m+1)(1-\alpha) \rceil$-largest value of $\{r_1,\cdots, r_m\}$.
\State \textbf{For the \emph{projected statistic}:}
\begin{align*}
 \stat(\theta) 
 = \frac{1}{n}\left(\pXX{n} - n\bp(\theta)\right)^\intercal  \Proj \hat M \Proj \left(\pXX{n} - n\bp(\theta)  \right)
\end{align*}
\State Set $\mintheta = \myargmin_{\theta \in \Theta} \stat(\theta)$
\State Sample $\{r_1,\cdots, r_m\}$ as  $m$ samples from the distribution of $\stat(\mintheta)$.
\State Set $\tau(\alpha,\epsilon)$ to be the $\lceil (m+1)(1-\alpha) \rceil$-largest value of $\{r_1,\cdots, r_m\}$.
\\ {\bf if } $\stat(\mintheta) > \tau(\alpha,\epsilon)$ {\bf then }
   Reject
\EndProcedure
\end{algorithmic}
\end{algorithm}

\subsection{Application - Goodness of Fit Testing}\label{sect:DP_GOF}
We first show that we can use the general chi-square test \texttt{DP-MC-MIN} with $\epsilon$-DP which uses Laplace noise in \Cref{alg:DP_MIN} for goodness of fit testing $H_0: \bp = \nullbp$.  In this case we select $\bp(\mintheta) = \nullbp$ and $\phi(\XX{n} + Z) = \nullbp$ in both the unprojected and projected statistics.  From the way that we have selected the critical value $\tau(\alpha,\epsilon)$ in  \Cref{alg:DP_MIN}, we have the following result on Type I error, which follows directly from Theorem 5.3 in \cite{GLRV16}.

\begin{theorem}
When the number of independent samples $m$ we choose for our MC sampling is larger than $1/\alpha$, testing $H_0:\bp = \nullbp$ in \Cref{alg:DP_MIN} guarantees Type I error at most $\alpha$ .
\end{theorem}

We then focus on empirically checking the power of \texttt{DP-MC-MIN} with $\alpha = 0.05$ for the different statistics.  As we did in the previous experiments, we will set the null hypothesis $\nullbp = (1/2,1/6,1/6,1/6)$ and alternate hypothesis $\bp^1 = \nullbp + 0.01 \cdot (1,-1/3,-1/3,-1/3)$ for various sample sizes.  We set the privacy parameter $\epsilon = \sqrt{2\cdot .001} \approx 0.045$, which implies $(\rho = 0.001)$-zCDP due to \Cref{thm:reduction}.  We set the number of independent samples we draw from the distribution of the statistic under the null hypothesis as $m = 59$.   In \Cref{fig:LAP_power_gof}, we compare the power of the projected and unprojected statistitic in $\texttt{DP-MC-GOF}$, labeled ``ProjGOF\_LAP" and ``UnProjGOF\_LAP" respectively, with the classical non-private chi-square test for various $n$ each with $5,000$ trials.  Note that there is a drastically larger power when we use the projected statistic as opposed to the unprojected statistic.  

We then show that the projected statistic using Laplace noise achieves significantly higher power than using the other unprojected test statistic as well as previous DP hypothesis tests with Laplace noise from \cite{GLRV16}.  We then label ``GLRV\_MCGOF\_LAP" as the MC based test with Laplace noise from \cite{GLRV16} and plot the power loss in \Cref{fig:LAPpower_gof_compare} that the other DP goodness of fit tests suffer when compared to the power that \texttt{DP-MC-GOF} achieves with the projected statistic.  Note that the error bars in the figure show $1.96$ times the standard error in the difference of proportions from $100,000$ trials.

\begin{figure}
\begin{center}
\begin{subfigure}{.45\textwidth}
\begin{flushleft}
\includegraphics[width=\linewidth]{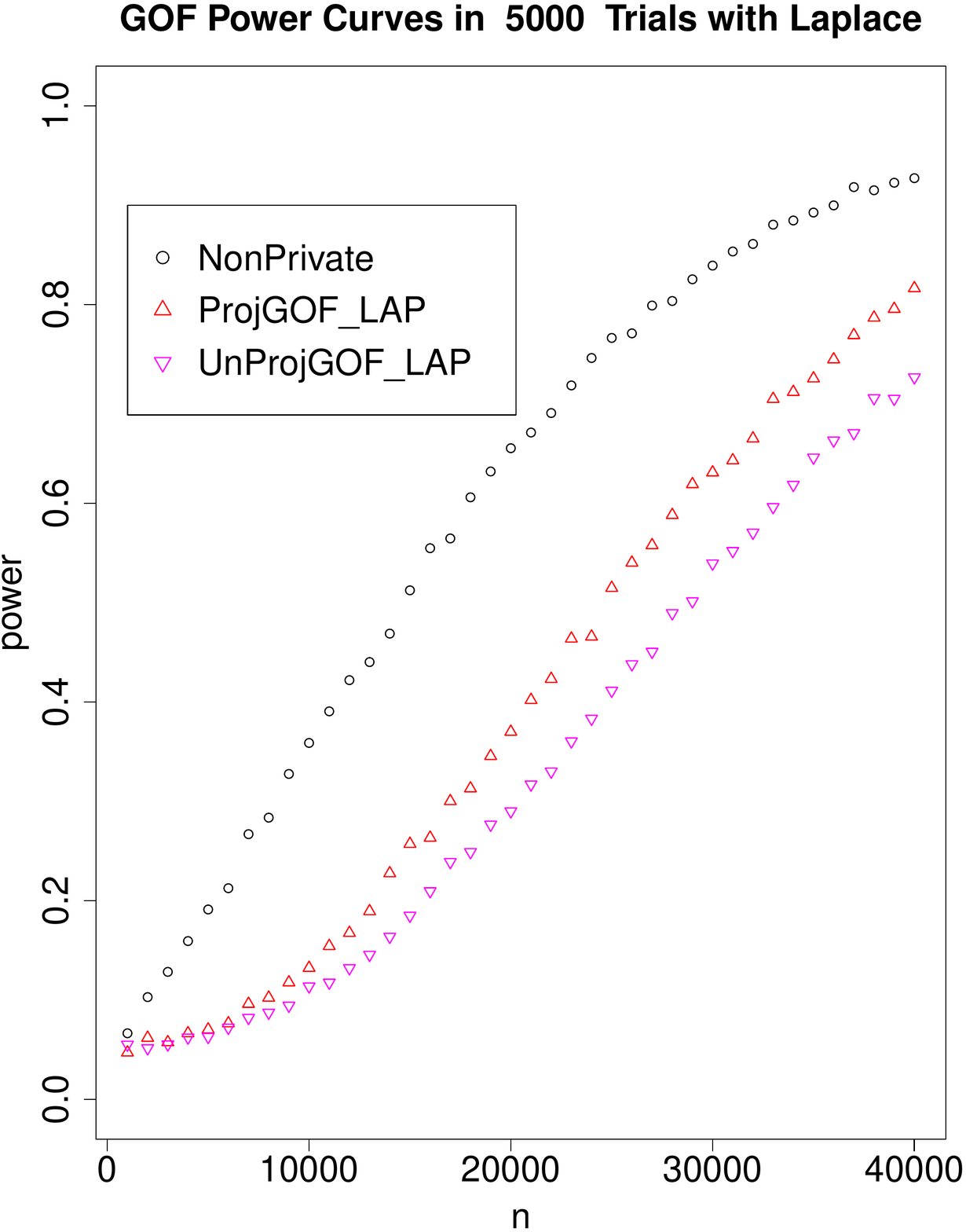}
\caption{A comparison of power for goodness of fit testing between the projected and unprojected statistics in $\texttt{DP-MC-MIN}$ with the classical nonprivate test for various $n$ with $5,000$ trials each. \label{fig:LAP_power_gof}}
\end{flushleft}
\end{subfigure}
\begin{subfigure}{.1\textwidth}
\end{subfigure}
\hspace{10mm}
\begin{subfigure}{.45\textwidth}
\begin{flushright}
\includegraphics[width=\linewidth]{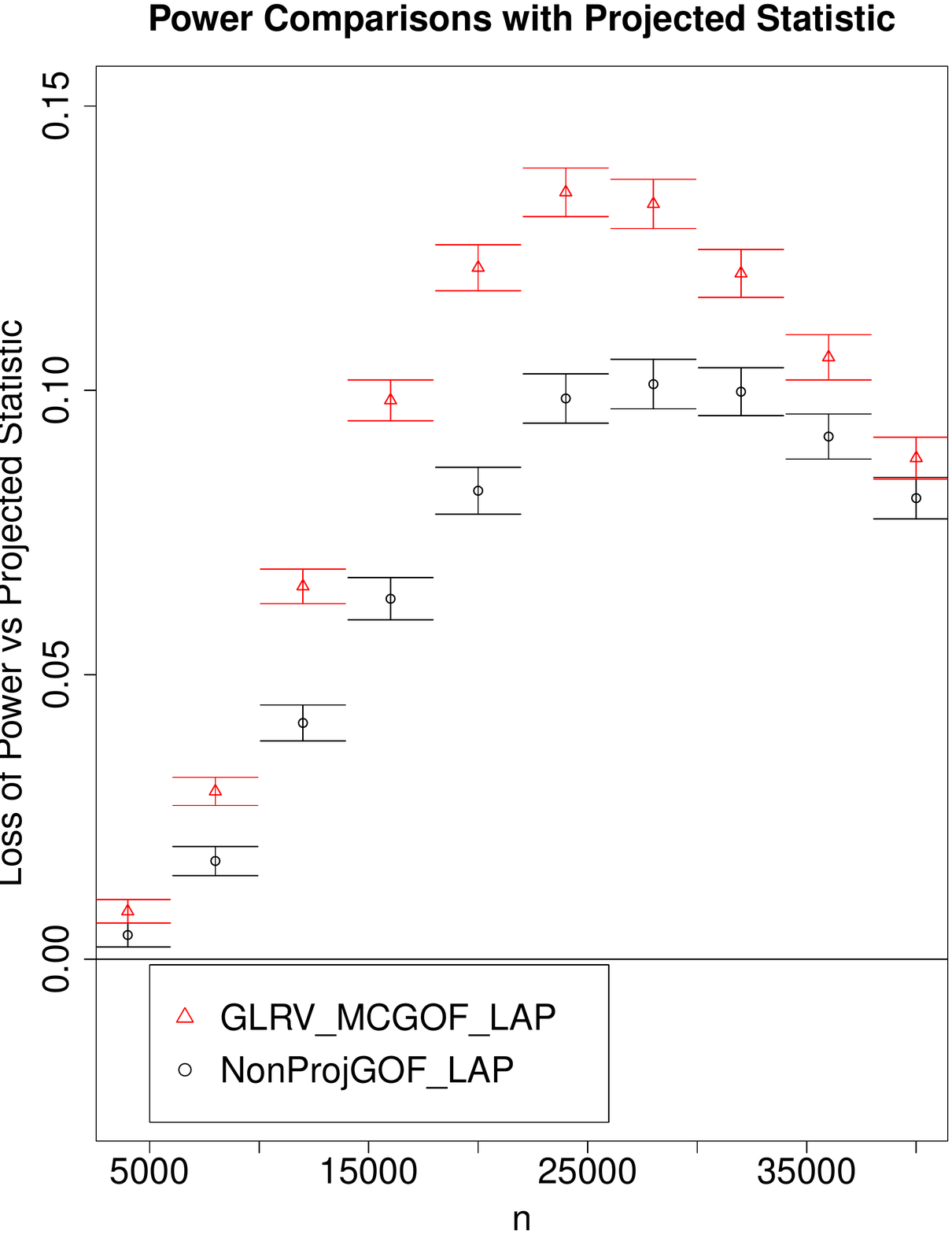}
\caption{The empirical loss in power from using other private goodness of fit tests instead of using the projected statistic in $\texttt{DP-MC-GOF}$.  The error bars correspond to $1.96$ times the standard error of each difference for $100,000$ trials. \label{fig:LAPpower_gof_compare}}
\end{flushright}
\end{subfigure}
\caption{Empirical power results for our new DP goodness of fit tests in $\texttt{DP-MC-MIN}$ and a comparisons to a previous private test in \cite{GLRV16} that uses Laplace noise with variance $4,000$ added to each cell count.}
\end{center}
\end{figure}

\subsection{Application - Independence Testing}
We then apply our general framework to independence testing as in \Cref{sect:IND}.  Unlike our goodness of fit testing, we are not guaranteed to have Type I error at most $\alpha$ when we have composite tests, e.g. independence testing, in $\texttt{DP-MC-MIN}$ because we are not sampling from the exact data distribution.   

We then empirically show the Type I error is at most the desired level $\alpha = 0.05$.  We again fix $\epsilon = \sqrt{2*0.001} \approx 0.045$, which ensures $\epsilon$-DP as well as $(\rho = 0.001)$-zCDP due to \Cref{thm:reduction}. We will use $m = 59$ samples in all of our MC testing.  We then give the empirical Type I error for various $n$ and data distributions in  \Cref{fig:DP_IND_signif}.  Note that we use the same rule of thumb as before where if our na\"ive estimate for the probability distribution produces expected cell counts smaller than $5$, then our test in inconclusive and fails to reject.  This is why in our experiments, the Type I error is close to zero for small sample sizes.  

\begin{figure}
\begin{center}
\begin{subfigure}{.45\textwidth}
\begin{flushleft}
\includegraphics[width=\linewidth]{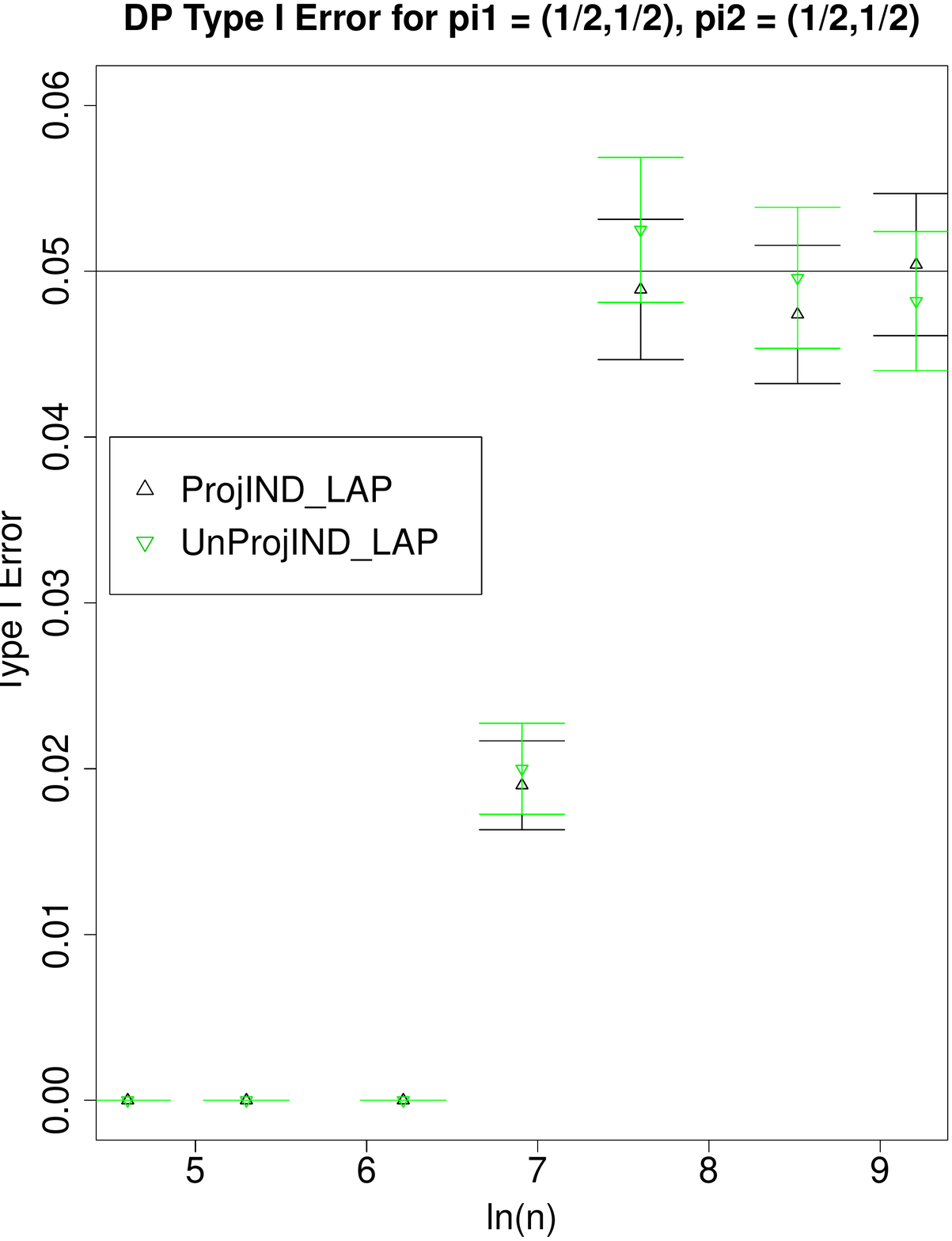}
\end{flushleft}
\end{subfigure}
\begin{subfigure}{.1\textwidth}
\end{subfigure}
\hspace{10mm}
\begin{subfigure}{.45\textwidth}
\begin{flushright}
\includegraphics[width=\linewidth]{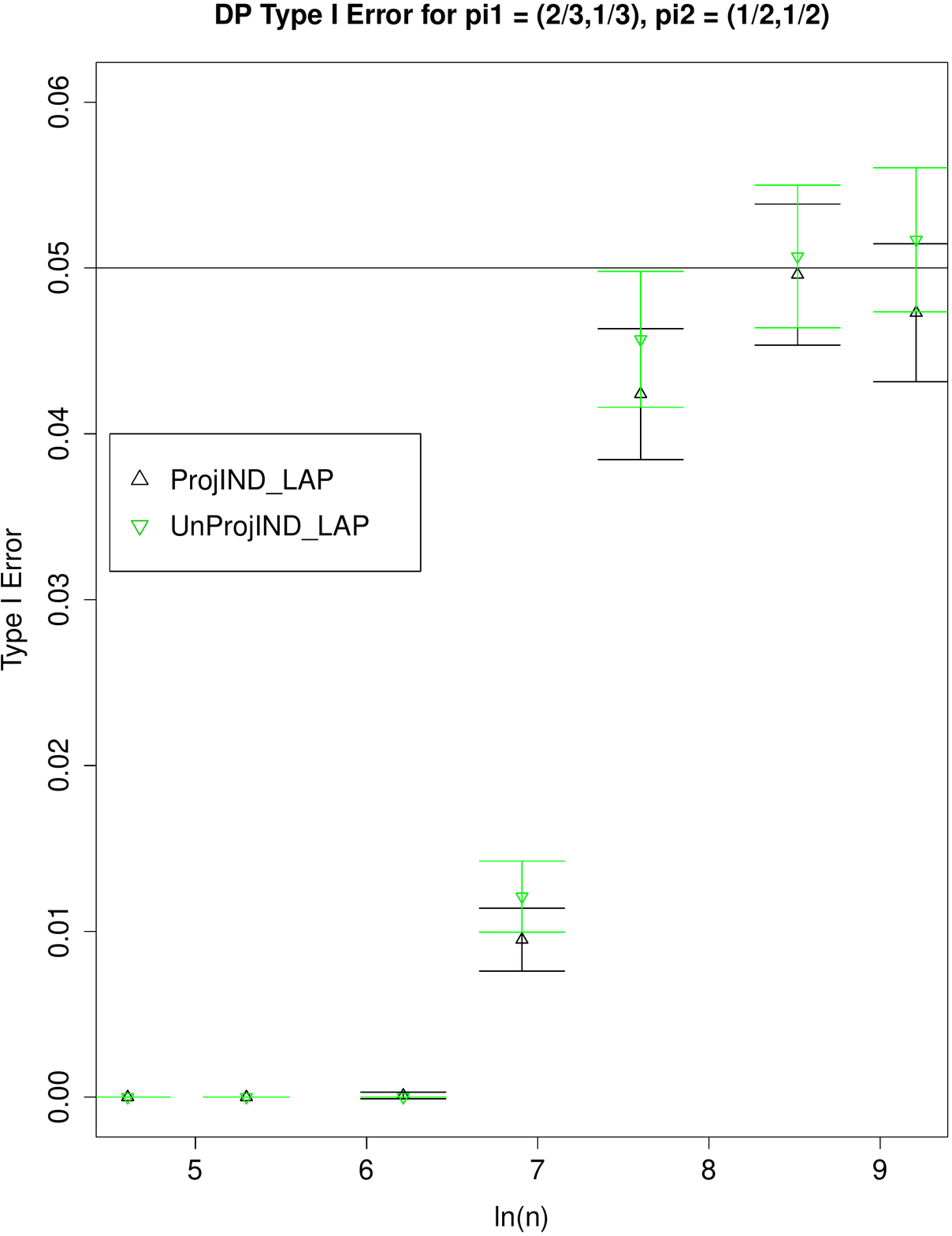}
\end{flushright}
\end{subfigure}
\caption{Empirical Type I Error for the new DP independence tests in $\texttt{DP-MC-MIN}$ with 1.96 times the standard error in $10,000$ trials.  We set $\epsilon = \sqrt{.001*2} \approx 0.045$ which corresponds to variance $1,000$ due to noise for each cell count.  It is desired to have Type I error at most $\alpha=0.05$, which is given as the horizontal line.  \label{fig:DP_IND_signif} }
\end{center}
\end{figure}

We also consider the power of our tests in $\texttt{DP-MC-MIN}$.   As we did before, we will set the data distribution with $\pipi{1} = (2/3,1/3), \pipi{2} = (1/2,1/2)$.  We then sample our contingency table $\XX{n}$ from $\mult(n,\bp(\pipi{1},\pipi{2})+ \pmb{\Delta})$ where $\pmb{\Delta} = 0.01 \cdot (1,0,-1,0)$ for various sample sizes.  In \Cref{fig:LAP_power_ind}, we compare the power of the projected and unprojected statistic in $\texttt{DP-MC-MIN}$, labeled ``ProjIND\_LAP" and ``UnProjIND\_LAP" respectively, with the classical non-private chi-square test for $1,000$ trials.  

We then label ``GLRV\_MCIND\_LAP" as the MC based test with Laplace noise from \cite{GLRV16}.  Note that the error bars show $1.96$ times the standard error in the difference of proportions from $10,000$ trials, giving a $95\%$ confidence interval.

\begin{figure}
\begin{center}
\begin{subfigure}{.45\textwidth}
\begin{flushleft}
\includegraphics[width=\linewidth]{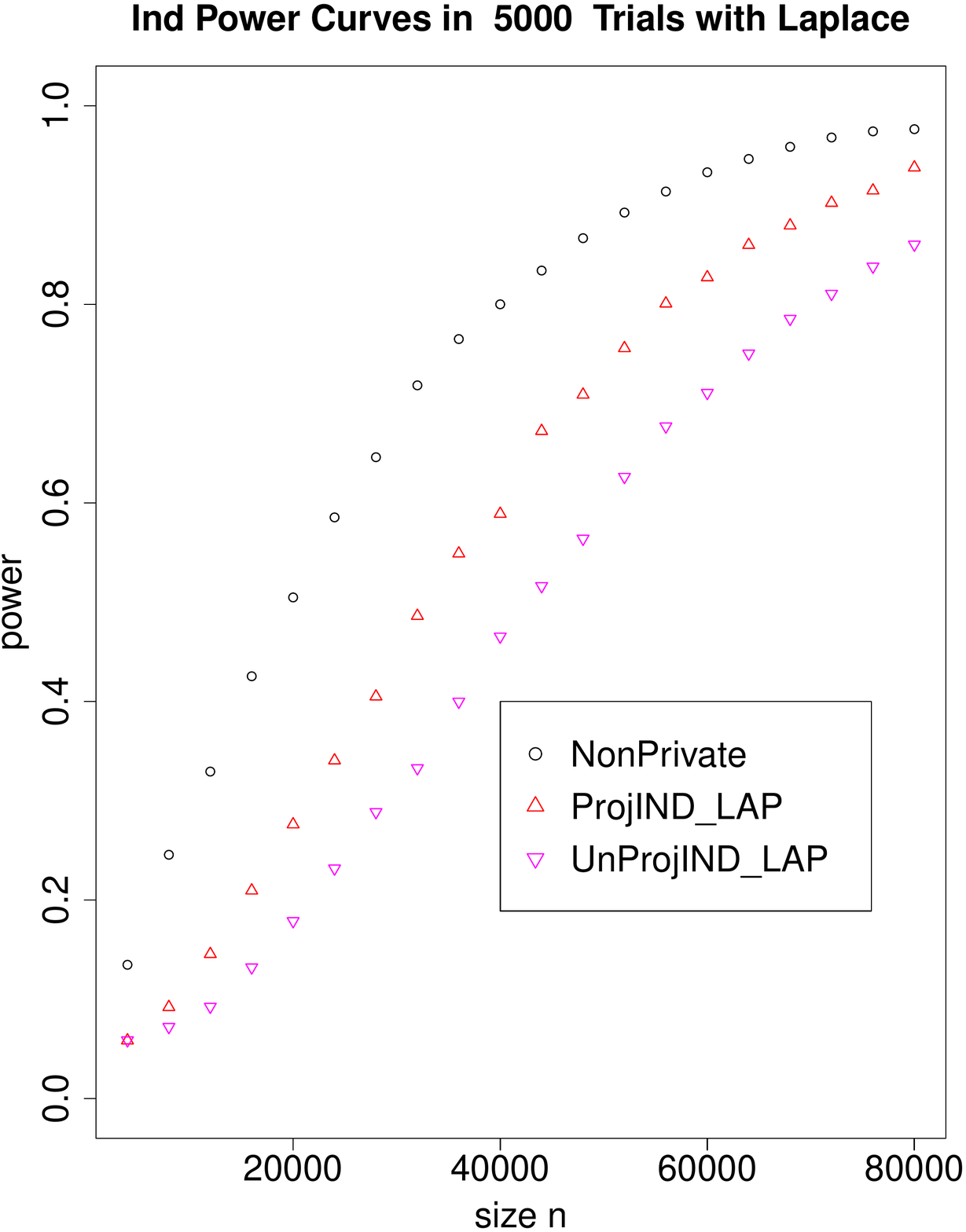}
\caption{A comparison of power for independence testing between the projected and unprojected statistics in $\texttt{DP-MC-MIN}$ with the classical nonprivate test.  \label{fig:LAP_power_ind}}
\end{flushleft}
\end{subfigure}
\hspace{10mm}
\begin{subfigure}{.1\textwidth}
\end{subfigure}
\begin{subfigure}{.45\textwidth}
\begin{flushright}
\includegraphics[width=\linewidth]{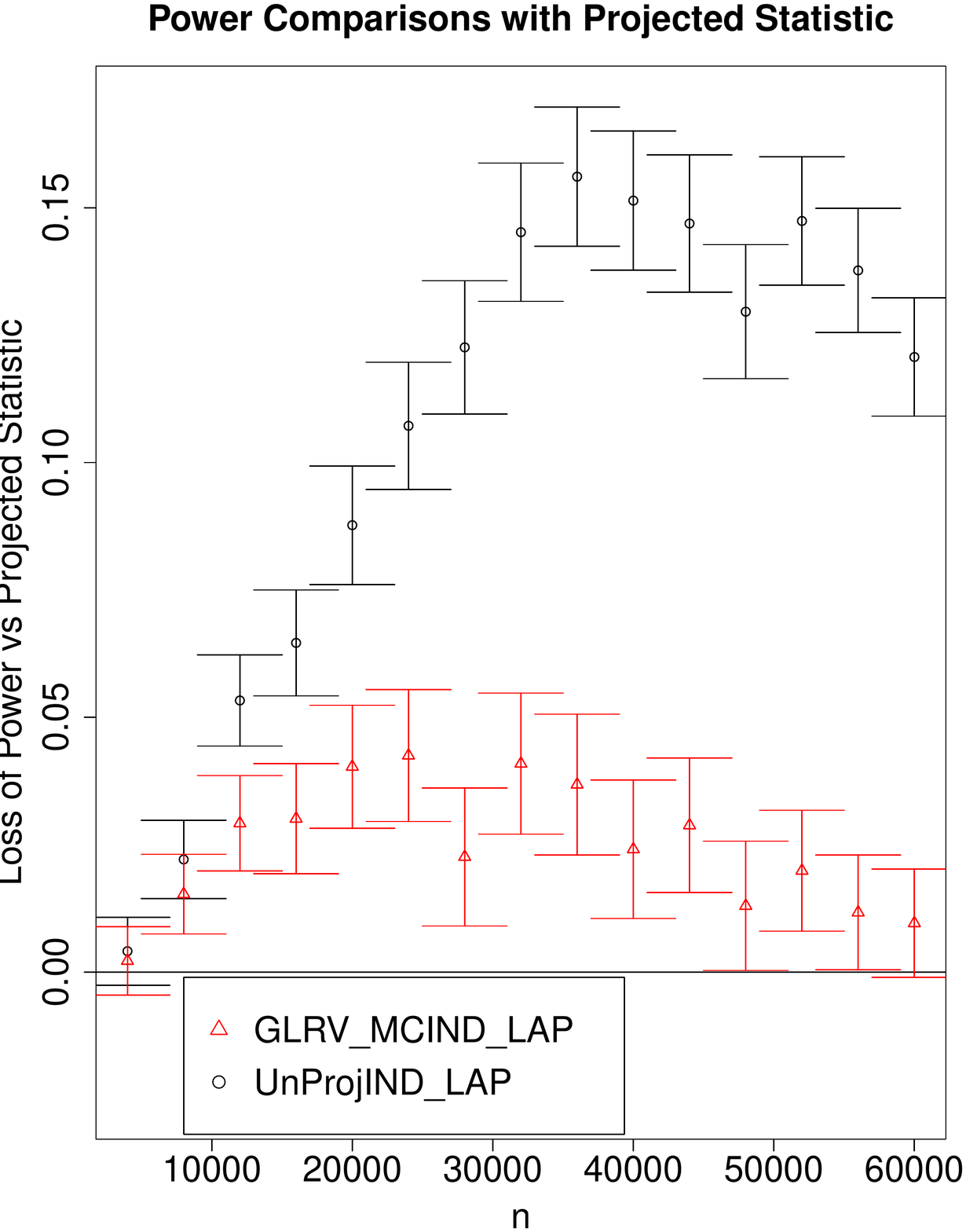}
\caption{The empirical loss in power from using other private goodness of fit tests instead of using the projected statistic in $\texttt{DP-MC-MIN}$ with error bars corresponding to $1.96$ times the standard error of each difference for $10,000$ trials. \label{fig:LAPpower_ind_compare}}
\end{flushright}
\end{subfigure}
\caption{Empirical power results for our new DP independence tests in $\texttt{DP-MC-MIN}$ and a comparisons to a previous private test in \cite{GLRV16} that uses Laplace noise.}
\end{center}
\end{figure}

\fi

\section{Conclusions}\label{sec:conclusions}
We have demonstrated a new broad class of private hypothesis tests $\MIN$ for categorical data based on the minimum chi-square theory.  We gave two statistics (\emph{unprojected} and \emph{projected}) that converge to a chi-square distribution when we use Gaussian noise and thus lead to zCDP hypothesis tests.  Unlike prior work, these statistics have the same asymptotic distributions in the private asymptotic regime as the classical chi-square tests have in the classical asymptotic regime.

Our simulations show that with either the unprojected or projected statistic our tests achieve at most $\alpha$ Type I error\ifnum\aistats=1 (see supplementary file)\fi.  We then empirically showed that our tests using the projected statistic significantly improves the Type II error when compared to the unprojected statistic and previous private hypothesis tests from \cite{GLRV16}.  Further, our new tests give comparable power to the classical (nonprivate) chi-square tests. \ifnum\aistats=1 See the supplementary file for \else We then gave \fi further applications of our new statistics to GWAS data and how we can incorporate other noise distributions (e.g. Laplace) using an MC sampling approach.


\clearpage
\bibliographystyle{abbrvnat}
\bibliography{./refs}
\clearpage
\begin{appendix}
\section{Proofs for \texorpdfstring{\Cref{subsec:minchi}}{Minimum Chi-Square Theory}}

\begin{proof}[Proof of \Cref{thm:min_chi}]
Since $\phi(\VV{n})$ converges in probability to $\truetheta$ and $M(\cdot)$ is a continuous mapping, then for any $b>0,c>0$ there exists an $n_0$ such that when $n\geq n_0$ then $M(\phi(\VV{n}))$ is within a distance $b$ from $M(\truetheta)$ with probability at least $1-c$, which makes $M(\phi(\VV{n}))$ positive definite with high probability for sufficiently large $n$.   
Furthermore, for any $d>0$, we can choose $n$ large enough so that the smallest eigenvalue of  $M(\phi(\VV{n}))$ is at least $\gamma-d$. 

Since the parameter space is compact, we know a minimizer exists for $\RR{n}(\theta)$.
Together, this implies that
 for sufficiently large $n$ and with high probability $\estQuadChi{n}(\mintheta)\geq 0$. 

Also, $\estQuadChi{n}(\mintheta)\leq \estQuadChi{n}(\truetheta)$ but $\estQuadChi{n}(\truetheta)/n\stackrel{P}{\rightarrow} 0$  since $M(\phi(\VV{n}))\stackrel{P}{\rightarrow} M$ and $\VV{n}\stackrel{P}{\rightarrow} A$.  Thus $\estQuadChi{n}(\mintheta)/n\stackrel{P}{\rightarrow} 0$ which means $\VV{n} - A(\mintheta)\stackrel{P}{\rightarrow} 0$ (since $M(\phi(\VV{n}))$ is positive definite with high probability and uniformly bounded away from $0$ in a neighborhood of $\truetheta$).  This implies that $A(\mintheta)\stackrel{P}{\rightarrow} A$  and so $\mintheta\stackrel{P}{\rightarrow} \truetheta$ since $A(\theta)$ is bicontinuous by assumption.

Thus, with high probability (e.g., $\geq 1-c$ for large enough $n$), $\mintheta$ satisfies the first order optimality condition $\nabla \estQuadChi{n}(\mintheta)=0$. This is the same as
\begin{align}
\dot{A}(\mintheta)^\intercal M(\phi(\VV{n})(\VV{n}-A(\mintheta))&=0\label{eqn:proof1}
\end{align}
Expanding $A(\mintheta)$ around $\truetheta$.
\begin{align}
A(\mintheta) &= A(\truetheta) + \underbrace{\int_0^1 \dot{A} (\truetheta + t(\mintheta - \truetheta))~dt}_{\equiv B(\mintheta)}(\mintheta - \truetheta)\label{eqn:proof2}
\end{align}
Substituting \eqref{eqn:proof2} into \eqref{eqn:proof1}, we get:
\begin{align}
& \dot{A} (\mintheta)^\intercal M(\phi(\VV{n}))\Big((\VV{n}-A(\truetheta)- B(\mintheta)(\mintheta-\truetheta)\Big)=0\\
& \dot{A} (\mintheta)^\intercal M(\phi(\VV{n}))B(\mintheta)\sqrt{n}(\mintheta-\truetheta)=\dot{A} (\mintheta)^\intercal M(\phi(\VV{n}))\sqrt{n}(\VV{n}-p(\truetheta) )
\end{align}
Now, by the continuity of $\dot{A}(\cdot)$ and the definition of $B$ and the convergence in probability of $\mintheta$ to $\truetheta$, we have $B(\mintheta)\stackrel{P}{\rightarrow} \dot{A} (\truetheta)$. Since $\dot{A}(\theta)$ has full rank by assumption, then for sufficiently large $n$,  $B(\mintheta)$ has full rank with high probability. This leads to the following expression with high probability for sufficiently large $n$,
\begin{align}
\sqrt{n}(\mintheta-\truetheta)&=\Big(\dot{A} (\mintheta)^\intercal M(\phi(\VV{n}))B(\mintheta)\Big)^{-1}\dot{A}(\mintheta)^\intercal M(\phi(\VV{n}))\sqrt{n}(\VV{n}-A)
\end{align}
Since $M(\phi(\VV{n}))$ has smallest eigenvalue at least $\gamma-d>0$ with high probability for $n$ large enough, and since $\phi(\VV{n})\stackrel{P}{\rightarrow} \truetheta$, $\mintheta\stackrel{P}{\rightarrow}\truetheta$, $B(\mintheta)\rightarrow \dot{A}(\truetheta)$ in probability, using  continuity in all of the above functions, and the assumption that $\sqrt{n}(\VV{n}-A)\rightarrow N(0,C)$ in distribution (and Slutsky's theorem) we get:
\begin{align}
\sqrt{n}(\mintheta - \truetheta)\stackrel{D}{\rightarrow} N\left(0, \Psi\right) \qquad \text{ as } \qquad n \to \infty.
\end{align}
\end{proof}

\begin{proof}[Proof of \Cref{thm:ferg}]
Note that Theorem 24 in \cite{Ferg96} shows that if the hypotheses hold then 
$$
n \left(\VV{n} - A(\mintheta) \right)^\intercal M(\mintheta ) \left( \VV{n} - A(\mintheta)\right) \stackrel{D}{\to} \chi^2_{\nu-k}.
$$
Note that we have $\phi(V^{(n)}) \stackrel{P}{\to} \truetheta$ and $\mintheta \stackrel{P}{\to} \truetheta$ for the true parameter $\truetheta \in \Theta$.  We can then apply Slutsky's Theorem due to $M(\cdot)$ being continuous, to obtain the result for $\estQuadChi{n}(\mintheta)$.
\end{proof}
\end{appendix}
\end{document}